\documentclass[12pt]{article}
\usepackage{amsmath}
\usepackage{amssymb}
\usepackage{amsthm}
\usepackage{amscd}
\usepackage{amsfonts}
\usepackage{graphicx}
\usepackage{fancyhdr}
\usepackage{color}
\usepackage{enumerate}
\usepackage{amsfonts}
\usepackage{amsthm}
\usepackage{psfrag} 
\usepackage{epsfig}
\usepackage{subfigure}
\usepackage{graphicx}
\usepackage{amssymb}
\usepackage{epstopdf}
\usepackage{amsbsy}
\usepackage{latexsym}
\usepackage{moreverb}                      
\usepackage{textcomp}

\usepackage{multicol}
\usepackage{enumitem}

\usepackage{algorithm}
\usepackage{algorithmic}

\newtheorem{remark}{Remark}

\newtheorem*{example*}{Example}

\DeclareMathOperator*{\argmin}{arg\,min}

\usepackage[top=2.8cm,bottom=2.8cm,left=2.5cm,right=2.5cm]{geometry}
\usepackage[T1]{fontenc}
\usepackage[utf8]{inputenc}

\usepackage[affil-it]{authblk} 
\usepackage{etoolbox}
\usepackage{lmodern}

\makeatletter
\patchcmd{\@maketitle}{\LARGE \@title}{\fontsize{18}{19.2}\selectfont\@title}{}{}
\makeatother

\usepackage{sectsty}
\sectionfont{\fontsize{13}{15}\selectfont}
\subsectionfont{\fontsize{12}{15}\selectfont}

\begin{document}

\title{A multi-fidelity neural network surrogate sampling method \\ for
  uncertainty quantification}

\author[1]{Mohammad Motamed\thanks{motamed@math.unm.edu}}

\affil[1]{Department of Mathematics and Statistics, The University of New Mexico, Albuquerque, USA}
\maketitle

\begin{abstract}
We propose a multi-fidelity neural network surrogate sampling method
for the uncertainty quantification of physical/biological systems
described by ordinary or partial differential equations. 
We first generate a set of low/high-fidelity data by low/high-fidelity
computational models, e.g. using coarser/finer discretizations of the
governing differential equations. 
We then construct a two-level neural network, where a large set
of low-fidelity data are utilized in order to accelerate the
construction of a high-fidelity surrogate model with a small set of
high-fidelity data. 
We then embed the constructed high-fidelity
surrogate model in the framework of Monte Carlo sampling. 
The proposed algorithm combines the approximation power of neural
networks with the advantages of Monte Carlo sampling within a multi-fidelity
framework. 
We present two numerical examples to demonstrate the accuracy and efficiency of the proposed
method. We show that dramatic savings in computational
cost may be achieved when the output predictions are desired to be accurate within small tolerances.

\end{abstract}

\medskip
\noindent
{\bf keywords:} 
multi-fidelity surrogate modeling, neural networks, uncertainty quantification


\section{Introduction}
\label{sec:intro}

Many physical and biological systems are mathematically modeled by
systems of ordinary and/or partial differential equations
(ODEs/PDEs). Examples include diffusion process, wave propagation, and DNA
transcription and translation, just to name a few. 
In addition to the possibility of involving multiple physical/biological
processes or multiple time/length scales, a major difficulty arises
form the presence of uncertainty in the systems and hence in the
ODE/PDE models. Uncertainty may be due to an inherent variability in
the system and/or our limited knowledge about the system, for instance
caused by noise and error in our measurements and experimental
data. The need for describing and quantifying uncertainty in complex
systems makes the field of uncertainty quantification (UQ) a
fundamental component of predictive science, enabling assertion of the
predictability of ODE/PDE models of complex systems; see e.g. \cite{UQ_book:15}. UQ problems are
often formulated in a probabilistic framework, where the models'
uncertain parameters are represented by probability distributions and
random processes. Solving the UQ problem then amounts to solving a
system of ODEs/PDEs with random parameters. 
In the present work we are particularly concerned with forward UQ
problems. Precisely, given
the distribution of uncertain parameters of the system, we want to
compute the statistics of some desired quantities of interest (QoIs).

There are a variety of methods available to tackle the forward UQ problem. 
An attractive class of methods, known as spectral methods (see
e.g. \cite{Ghanem,NTW:08,Xiu_Hesthaven}), exploits the possible
regularity that output QoIs might have with respect to the input
parameters. The performance of these methods, however, dramatically
deteriorates in the absence of regularity. For example, solutions to
parametric hyperbolic PDEs are in general non-smooth, and therefore
related stochastic QoIs are often not regular; see
e.g. \cite{Motamed_etal:13,Motamed_etal:15}. Consequently, spectral
methods may not be applicable to stochastic hyperbolic problems. 
Another popular method for solving the UQ problem is Monte Carlo (MC) sampling \cite{MC}. While
being flexible and easy to implement, MC sampling features a very
slow convergence rate. More recently, a series of advanced variants of MC sampling has been proposed to speed up computations; see
e.g. \cite{Giles:08,Giles:11,MOMC:18,Multiindex:16,QMC:11,MLQMC:15,Speight:09,CVMLMC_Nobile:15,Gorodetsky_etal:19}
and the references there in. 
%
%
%

In recent years, neural networks have shown remarkable success in
solving a variety of large-scale artificial intelligence problems; see
e.g. \cite{Schmidhuber:15}. 
Importantly, the application of neural networks is not limited to artificial
intelligence problems. Neural networks can, for instance, be utilized to build surrogate models
for physical/biological QoIs and hence to solve UQ problems. 
In this setting, the data to construct (or {\it train}) a neural
network as a surrogate model for a QoI needs to be obtained by computing realizations
of the QoI, where each realization involves solving an ODE/PDE problem through numerical discretization. 
We may however face a problem with this approach: the construction of an accurate
neural network surrogate model usually requires abundant high-fidelity
data that in turn amounts to a large number of high-fidelity computations that
may be very expensive or even infeasible.

Motivated by multi-fidelity approaches (see e.g. \cite{MF_review:16}), 
we propose to construct a two-level neural network, 
where a large set of low-fidelity data are utilized in order to
accelerate the construction of a high-fidelity surrogate model with a
small set of high-fidelity data. We then embed the constructed
high-fidelity surrogate model in the framework of MC sampling for
solving the UQ problem in hand. 
The main goal is to combine the approximation power of neural
networks with the advantages of MC sampling within a multi-fidelity
framework. 
More precisely, we construct an algorithm consisting of the following steps. 
We first generate a large set of low-fidelity and a
small set of high-fidelity data by low/high-fidelity computational
models, e.g. using coarser/finer discretizations of the governing
ODEs/PDEs. 
We then construct two networks at two different levels. 
The network in the first level, denoted by NN1,
uses all available data to capture the (possibly nonlinear) {\it
  correlation} between the two computational models at different levels of fidelity. 
The trained network NN1 will then be utilized to produce additional high-fidelity data at
a very small cost. 
Next, the network in the second level, denoted by NN2, uses all
available and newly generated high-fidelity data to construct a
surrogate model as an accurate prediction of the QoI. 
Finally, the statistical moments of the QoI are approximated by its sample
moments, where each realization is
computed through a fast evaluation of
the constructed high-fidelity surrogate model. 
We present two numerical examples to motivate the efficiency of the
proposed algorithm compared to MC sampling. In particular, we show
that we may achieve dramatic savings in computational cost when the output predictions are desired to be accurate within small tolerances.

It is to be noted that the framework presented here is not limited to two
levels of fidelity and can hence be extended to training more than two
networks using data sets at multiple levels of fidelity. 
The construction can
also be adapted to the more advanced variants of MC sampling, such as
multi-level, multi-order, and multi-index methods (see
e.g. \cite{Giles:11,CVMLMC_Nobile:15,MOMC:18,Multiindex:16}) and to
other sampling-based techniques such as multi-index collocation (see
e.g. \cite{Adaptive_MISC:19}). 
Moreover, the proposed sampling method may also be applied to inverse UQ problems for the
efficient computation of marginal likelihoods. 
We refer to \cite{Aydin_etal:19, Liu_Wang:19,MFNN:20} for other recent
works in the context of multi-fidelity neural networks; see Remark
\ref{Composite_MFNN} for a short comparison between the proposed
construction here and the construction in \cite{MFNN:20}.

The rest of the paper is organized as follows. In Section \ref{sec:uq}
we state the mathematical formulation of the UQ problem and briefly
address its numerical treatment. 
We then present the multi-fidelity neural network surrogate
sampling algorithm in Section \ref{sec:mfnn} and discuss its cost and accuracy. 
In Section \ref{sec:numerics} we perform two numerical examples: and ODE problem and a
PDE problem. Finally, in Section \ref{sec:conclusion} we summarize conclusions and outline future works.

\section{Uncertainty Quantification of Physical Systems}
\label{sec:uq}

In this section, we first present the mathematical formulation of the
UQ problem that we consider in the present work. We then briefly
review the available numerical methods.

\subsection{Problem statement}

Let $M_{DE}$ denote a mathematical model consisting of a set of
parametric ODEs/PDEs with a $p$-dimensional uncertain parameter vector $\Theta \in {\mathbb R}^p$ that is represented by an $n$-dimensional
random vector ${\bf y} \in \Gamma \subset {\mathbb R}^{n}$ with a known and bounded joint PDF $\pi:
\Gamma \rightarrow {\mathbb R}_+$. 
Suppose that we want to map the random input parameters through the
model to obtain a desired output quantity $Q \in {\mathbb R}$. 
In abstract form we can write 
$$
Q({\bf y}) = M_{DE}(\Theta({\bf y})).
$$
%
%
Since $Q$ is a random quantity in the
light of randomness in $\Theta$, our specific goal is to compute the
statistics of $Q$. For instance, we may be interested in computing the first statistical
moment (or expectation) of $Q$,
$$
{\mathbb E}[Q] := \int_{\Gamma} Q({\bf y}) \, \pi({\bf y}) \, d{\bf y}.
$$

We remark that this may be a challenging problem especially when the
quantity $Q$ is not highly regular with respect to ${\bf y}$ or
when ${\bf y}$ lives in a high-dimensional space, i.e. when $n$ is
large. 
In such cases, an
accurate approximation of ${\mathbb E}[Q]$ may require many evaluations of
$M_{DE} (\Theta)$ corresponding to many realizations of $\Theta$, where each
evaluation involves computing a complex (and often expensive) DE problem. 
%


\subsection{Numerical methods} 
\label{sec:methods}

A popular method for computing the statistics of $Q$ is MC sampling \cite{MC}, where sample statistics of $Q$ are computed from
independent realizations drawn from the distribution $\pi$. Two
favorable features of MC sampling include its flexibility with respect
to irregularity and high dimensionality. MC sampling can easily handle situations
where the map $Q: \Gamma \rightarrow {\mathbb R}$ is not regular and
$\Gamma \subset {\mathbb R}^n$ is a high-dimensional space. Despite
these advantages, MC smapling features a
very slow convergence rate. This may make MC simulations infeasible,
since a very large number of (possibly expensive) DE problems are
required to be solved in order to obtain accurate approximations. 

Another class of methods, such as stochastic Galerkin \cite{Ghanem}
and stochastic collocation \cite{NTW:08,Xiu_Hesthaven}, 
can exploit the possible regularity that the map $Q: \Gamma \rightarrow {\mathbb R}$ might have. These methods are expected to
yield a very fast spectral convergence provided $Q$ is highly regular. The performance of spectral methods, however, strongly deteriorates in the
presence of low regularity and high dimensionality; see
e.g. \cite{Motamed_etal:13,Motamed_etal:15}. 
%
%

More recently, several variants of MC sampling have been proposed that
retain the two aforementioned advantages of MC sampling and yet accelerate its slow convergence. Examples include
multi-level MC \cite{Giles:08,Giles:11}, multi-order MC \cite{MOMC:18}, multi-index MC
\cite{Multiindex:16}, quasi MC
\cite{QMC:11}, multi-level quasi MC \cite{MLQMC:15}, control variate
multi-level MC \cite{Speight:09,CVMLMC_Nobile:15}, and control variate multi-fidelity
MC \cite{Gorodetsky_etal:19}, just to name a few. It is to be noted
that while multi-index and quasi MC approaches require mild
assumptions on the regularity of $Q$, the rest of the methods in this
category are resilient with respect to regularity.

In what follows we present a new approach that combines the
approximation power of neural networks with the advantages of MC
sampling in a multi-fidelity framework. This method may be
considered as an advanced varient of MC sampling and classified in the
third category listed above.

\section{Multi-fidelity Neural Network Surrogate Sampling}
\label{sec:mfnn}

In this section we will first review the main notions of
multi-fidelity modeling relevant to the focus of this
work 
and give a brief overview of feedforward neural networks used in
surrogate modeling. 
%
We then present a multi-fidelity neural network
surrogate sampling algorithm for uncertainty quantification and
discuss its accuracy and computational cost.

\subsection{Selection of multi-fidelity models and their correlation}

Let $Q_{LF}({\bf y})$ and $Q_{HF}({\bf y})$ denote the
approximated values of the quantity $Q({\bf y})$ by a low-fidelity and a
high-fidelity computational model, respectively. We make the following
assumptions on the two computational models:

\begin{itemize}
\item[A1.] The high-fidelity quantity $Q_{HF}({\bf y})$ is an accurate approximation of
the quantity $Q({\bf y})$ within a desired small tolerance.

\item[A2.] The low-fidelity quantity $Q_{LF}({\bf y})$ is another
  approximation of $Q({\bf y})$ that is both correlated with $Q_{HF}({\bf
    y})$ and computationally
  cheaper than computing $Q_{HF}({\bf y})$.
\end{itemize}

The high-fidelity model is often obtained by a direct and fine
discretization of the underlying model $M_{DE}$. 
There are however several possibilities to build the low-fidelity
model. For example, we may build a low-fidelity model by directly
solving the original ODE/PDE problem using either a coarse discretization or a low-rank approximation of
linear systems that appear in the high-fidelity model. Another
possibility for building a low-fidelity model is to solve an auxiliary
or effective problem obtained by simplifying the original ODE/PDE
problem. For instance we may consider a model with simpler
physics, or an effective model obtained by homogenization, or a
simpler model obtained by smoothing out the rough parameters of
$M_{DE}$ or through linearization. 
Without the loss of generality, here we consider a coarse discretization and
a fine discretization of the underlying model $M_{DE}$ to build up
the low-fidelity and high-fidelity computational models,
respectively. 
We hence write
\begin{equation}\label{bi-fidelity}
Q_{LF}({\bf y}) = M_{DE}^{h_{LF}}(\Theta({\bf y})), \qquad Q_{HF}({\bf
  y}) = M_{DE}^{h_{HF}}(\Theta({\bf y})), \qquad h_{HF} < h_{LF},
\end{equation}
where $h_{LF}$ and $h_{HF}$ denote the mesh size and/or the time step
of a stable discretization scheme used in the low-fidelity and high-fidelity models, respectively. We also assume that the
error in the approximation of $Q$ by the high-fidelity model satisfies
\begin{equation}\label{HF_error}
| Q({\bf y}) - Q_{HF}({\bf y}) | = |M_{DE}(\Theta({\bf y})) -
M_{DE}^{h_{HF}}(\Theta({\bf y})) | \le c \, h_{HF}^q, \qquad 
\forall {\bf y} \in \Gamma,
\end{equation}
where $q>0$ is related to the order of accuracy of the discretization
scheme, and $c>0$ is a bounded constant. This upper error bound implies that A1 holds. Moreover, for A2 to hold, we further assume
that $h_{LF}$ (while being larger than $h_{HF}$) 
is small enough for $Q_{LF}$ to be correlated with $Q_{HF}$. 

It is to be noted that in general the selection of low-fidelity and high-fidelity models
is problem dependent. A bi-fidelity model is admissible as long as
assumptions A1-A2 are satisfied, although it may not be the best
bi-fidelity model in terms of computational efficiency while
achieving a desired accuracy among all admissible bi-fidelity models. 


A major step in multi-fidelity modeling is to capture and utilize the correlation
between the models at different levels of fidelity. A widely used
approach, known as comprehensive
correction, is to assume a linear correlation and write $Q_{HF}({\bf y}) =
\rho({\bf y}) Q_{LF}({\bf y}) + \delta({\bf y})$, where $\rho$ and $\delta$ are the unknown
multiplicative and additive corrections, respectively; see
e.g. \cite{MF_review:16}. The main limitation of this
strategy is its inability to capture a possibly nonlinear relation
between the two models. Therefore, we consider a general correlation and express the relation between the two
models as
\begin{equation}\label{corr_function}
Q_{HF}({\bf y}) = F({\bf y}, Q_{LF}({\bf y})),
\end{equation}
where $F$ is a general unknown function that captures the (possibly
nonlinear) relation between the low-fidelity and high-fidelity quantities; see also \cite{MF_nonlinear:17,MFNN:19}. 
%

\subsection{Feedforward neural network surrogate modeling}
\label{sec:NN}

\begin{figure}[!h]
\begin{center}
\includegraphics[width=0.45\linewidth]{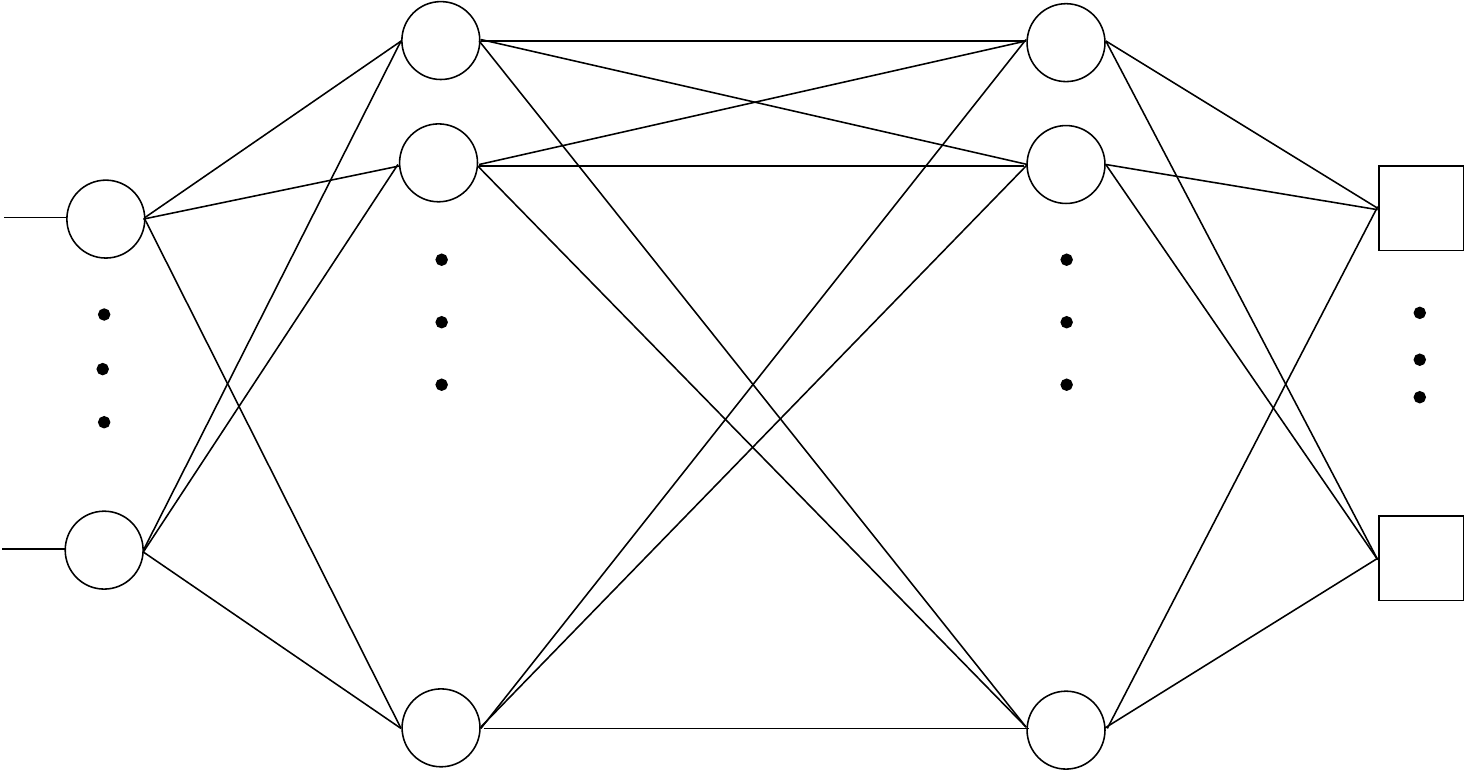}     
\caption{Graph representation of a feed-forward network with $L=2$
  hidden layers.}
\label{graph_fig}
\end{center}
\end{figure}

Consider a feedforward neural network with one input layer consisting
of $n_{\text{in}} \in {\mathbb N}$ neurons, one output layer consisting
of $n_{\text{out}} \in {\mathbb N}$ neurons, and $L \in
{\mathbb N}$ hidden layers consisting of $n_{1}, \dotsc, n_{L} \in {\mathbb
  N}$ neurons, respectively. A neural network with such an
architecture may be represented
by a map ${\bf f}_{\text{NN}}: {\mathbb
  R}^{n_{\text{in}}} \rightarrow {\mathbb R}^{n_{\text{out}}}$ given
by the composition 
$$
{\bf f}_{\text{NN}}(\boldsymbol\theta) = {\bf f}_{L+1} \circ {\bf f}_L
\circ \cdots \circ {\bf f}_1(\boldsymbol\theta), \qquad
\boldsymbol\theta \in {\mathbb R}^{n_{\text{in}}},
$$ 
where each individual map ${\bf f}_{\ell}: {\mathbb R}^{n_{\ell - 1}} \rightarrow
{\mathbb R}^{n_{\ell}}$, with $\ell = 1, \dotsc, L+1$, $n_0 =
n_{\text{in}}$, and $n_{L+1} =
n_{\text{out}}$, is obtained by the component-wise application of
a (possibly nonlinear) activation function $\sigma_{\ell}$ to an affine-linear map 
$$
{\bf f}_{\ell}({\bf z}) = \sigma_{\ell} (W_{\ell} \, {\bf z} + {\bf
  b}_{\ell}), \qquad {\bf z} \in {\mathbb R}^{n_{\ell-1}}, \quad 
W_{\ell} \in {\mathbb R}^{n_{\ell} \times n_{\ell-1}}, \quad
{\bf b}_{\ell} \in {\mathbb R}^{n_{\ell}}, \quad 
\ell = 1, \dotsc, L+1.
$$
The parameters $W_{\ell}$ and ${\bf b}_{\ell}$ that define the
affine-linear maps in layer $\ell$ are referred to as the weights (or
edge weights) and
biases (or node weights) of the $\ell$-th layer, respectively. 
Popular choices of activation function include the hyperbolic tangent
$\sigma_{\ell}(\theta) = \tanh (\theta)$, the sigmoid function
$\sigma_{\ell}(\theta) = 1/ (1+ \exp(- \theta))$, and the rectified linear
unit (ReLU) activation function $\sigma_{\ell} (\theta) = \max (0, \theta)$, where $\theta \in
{\mathbb R}$. Different layers of a network may use different
activation functions. In surrogate modeling applications, the identity
function is often used as the activation function in the output
layer, i.e. $\sigma_{L+1}(\theta)=\theta$. 
Figure \ref{graph_fig} shows a
graph representation of the network, where each node represents a
neuron, and each edge connecting two nodes represents a multiplication by a
scalar weight. The input neurons/nodes takes the $n_{\text{in}}$ components of
the independent variable $\boldsymbol\theta=(\theta_1, \dotsc,
\theta_{n_{\text{in}}})$, and the output neurons/nodes produce the
$n_{\text{out}}$ components of ${\bf f}_{NN} = (f_{NN,1}, \dotsc, f_{NN,n_{\text{out}}})$.

In the context of surrogate modeling, we aim at {\it training} a
network ${\bf f}_{\text{NN}}$ with
pre-assigned architecture and activation functions that {\it learns}
an unknown function ${\bf f}: {\mathbb
  R}^{n_{\text{in}}} \rightarrow {\mathbb R}^{n_{\text{out}}}$. 
Precisely, we want to find a set of network parameters $\Theta_{\text{NN}} := \{ (W_{\ell}, {\bf
  b}_{\ell}) \}_{\ell=1}^{L+1}$ such that ${\bf f}_{\text{NN}}
(\boldsymbol\theta; \Theta_{\text{NN}})$ well approximates ${\bf f}(\boldsymbol\theta)$ within a small
tolerance. Here, by an abuse of notation, we write ${\bf f}_{\text{NN}}(\boldsymbol\theta;
\Theta_{\text{NN}})$ to emphasize the dependence of the network on
the parameter set $\Theta_{\text{NN}}$. 
Such a trained network may then serve as a
surrogate model for the desired function ${\bf f}$ that may
be only indirectly available through, for instance, a set of complex ODE/PDE problems. 

For this purpose, we first collect a data
set of $M$ input-output pairs $\{
(\boldsymbol\theta^{(i)}, {\bf f}(\boldsymbol\theta^{(i)})) \}_{i=1}^M$ and
then formulate an optimization problem
$$
\argmin_{\Theta_{\text{NN}}} \frac{1}{M}
\sum_{i=1}^M C({\bf f}_{\text{NN}}
(\boldsymbol\theta^{(i)}; \Theta_{\text{NN}}), {\bf f}(\boldsymbol\theta^{(i)})),
$$
where $C({\bf a}, {\bf b})$ is a cost (or {\it loss}) function that measures the distance
between vectors ${\bf a}$ and ${\bf b}$. For instance, a typical cost function is the quadratic cost
function $C({\bf a}, {\bf b}) = || {\bf a} - {\bf b} ||^2$. 
The optimization problem may then be solved by a gradient-based method, such as stochastic gradient
descent \cite{SGD1,SGD2} or Adam \cite{Adam:17}. In these methods the
gradient of the cost function with
respect to the network parameters is usually computed by the
chain rule using a differentiation technique known as {\it back
  propagation} \cite{BackProp:1986}. We refer
to the review paper \cite{Bottou_etal:18} for more
details. 

It is to be noted that the above training strategy may suffer from
{\it overfitting}, that is, the trained network may not perform well
in approximating ${\bf f}(\boldsymbol\theta)$ outside the set of input training data $\{\boldsymbol\theta^{(i)}\}_{i=1}^M$. In order to avoid
  overfitting, a few {\it regularization}
  techniques have been proposed. Examples include the addition of a
  regularization (or penalty) term to the cost function, early stopping, and
  dropout \cite{Dropout:14}. 
Each optimization/regularization technique involves a few
hyper-parameters, such as learning rates, number of epochs, mini-batch
sizes, penalty parameters, dropout
rates, and so forth. The hyper-parameters are often tuned using a {\it validation}
set, i.e. a set of data points that are not directly used in the
optimization process. A common practice is to select a large portion
of the available set of $M$
data points as training data and a smaller portion of the data as the validation
data. We refer to \cite{Bengio:12,GoodBengCour16} for more details on
the subject.

The performance of a surrogate neural network, i.e. the accuracy of
the approximation ${\bf f}(\boldsymbol\theta) \approx {\bf f}_{\text{NN}}
(\boldsymbol\theta; \Theta_{\text{NN}})$, would depend on two factors:
1) the structure of the network, that is, the number of layers and neurons,
activation functions, and the optimization/regularization techniques
and parameters used to train the network, and 
2) the choice of the data set $\{(\boldsymbol\theta^{(i)}, {\bf f}(\boldsymbol\theta^{(i)})) \}_{i=1}^M$. 
Currently, there is no rigorous theory addressing the precise
dependence of the accuracy of the approximation on the
aforementioned factors. Much of the work has been focused on the {\it
expressivity} of neural networks, i.e. the ability of neural networks
in terms of their structure to approximate a
wide class of functions. These works range from the universal
approximation theorem for feed-forward networks with one hidden layer
(see e.g. \cite{Hornik_etal:89}) to more recent results for networks
with more complex structures and wider classes of functions; see
e.g. \cite{MhaskarPoggio:16,Yarotsky:17,PetVio:18,MontanelliDu:19,SchwabZech:19,Bolcskei_etal:19,Guhring_etal:19}
and the references therein. 
While showing the approximation power of
neural networks, these results do not enable a cost-error analysis to
achieve a desired accuracy with minimum computational cost. 
%
%
Moreover, when the target function ${\bf f}$ is indirectly given by a set of complex
ODE/PDE problems, we need to compute an accurate approximation of the function at the input data points
$\{ \boldsymbol\theta^{(i)} \}_{i=1}^M$ to obtain the output data
points $\{ {\bf f}(\boldsymbol\theta^{(i)}) \}_{i=1}^M$. In such cases, the
choice of the data set, i.e. the number $M$ and location of the input points $\{
\boldsymbol\theta^{(i)} \}_{i=1}^M$ in the domain space and the accuracy of the output
data points $\{ {\bf f}(\boldsymbol\theta^{(i)}) \}_{i=1}^M$, is crucial in
establishing an accurate and efficient neural network surrogate
model. An important and open question that arises is then: given a
fixed structure, how well is the approximation with respect to the
choice of data points?


\subsection{A multi-fidelity neural network surrogate-based sampling algorithm}
\label{sec:algorithm}

We first generate a large set of low-fidelity and a small set of high-fidelity data by low/high-fidelity computational models. We then construct two networks at two different levels. 
The first neural network, denoted by NN1, uses all available data to perform two
tasks. First, it learns the correlation between
low-fidelity and high-fidelity data. Second, the learned correlation
is used to generate extra high-fidelity data. The second neural
network, denoted by NN2, is then trained using the original and newly generated
high-fidelity data to serve as a surrogate for the high-fidelity
quantity. 
The constructed surrogate model is then embedded within the framework of MC sampling to
compute the statistics of the QoI. 

\begin{figure}[!h]
\begin{center}
\includegraphics[width=0.65\linewidth]{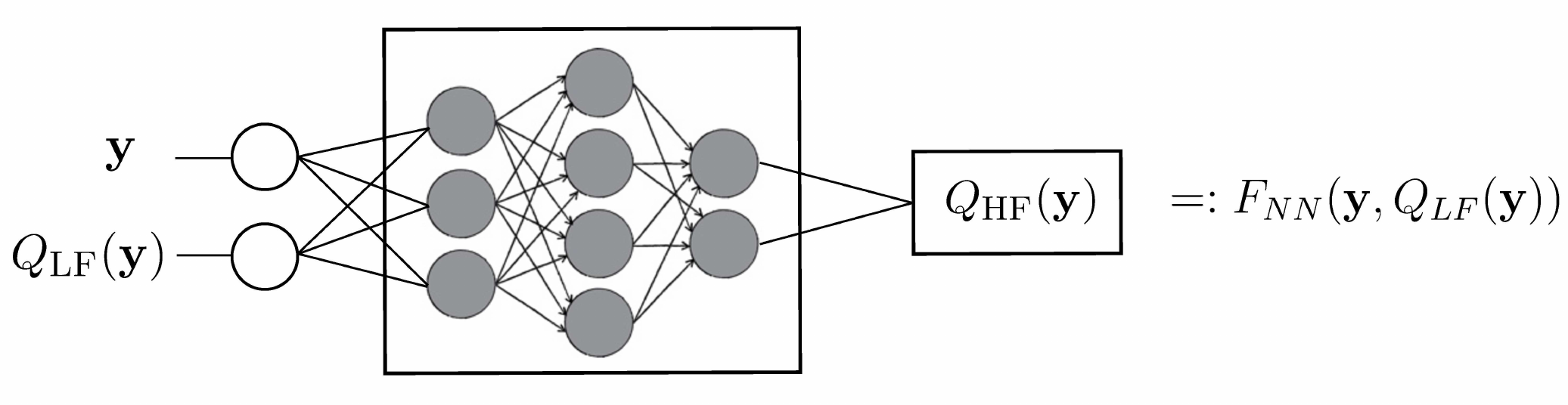}     
\caption{A schematic representation of NN1. It takes the input data
  $\{ ( {\bf y}^{(i)}, Q_{LF}^{(i)} )\}_{i=1}^{M_1}$ and the output data
  $\{Q_{HF}^{(i)} \}_{i=1}^{M_1}$ to learn the correlation function $F({\bf y},
  Q_{LF}({\bf y}))$.}
\label{NN1_fig}
\end{center}
\end{figure}

\begin{figure}[!h]
\begin{center}
\includegraphics[width=0.53\linewidth]{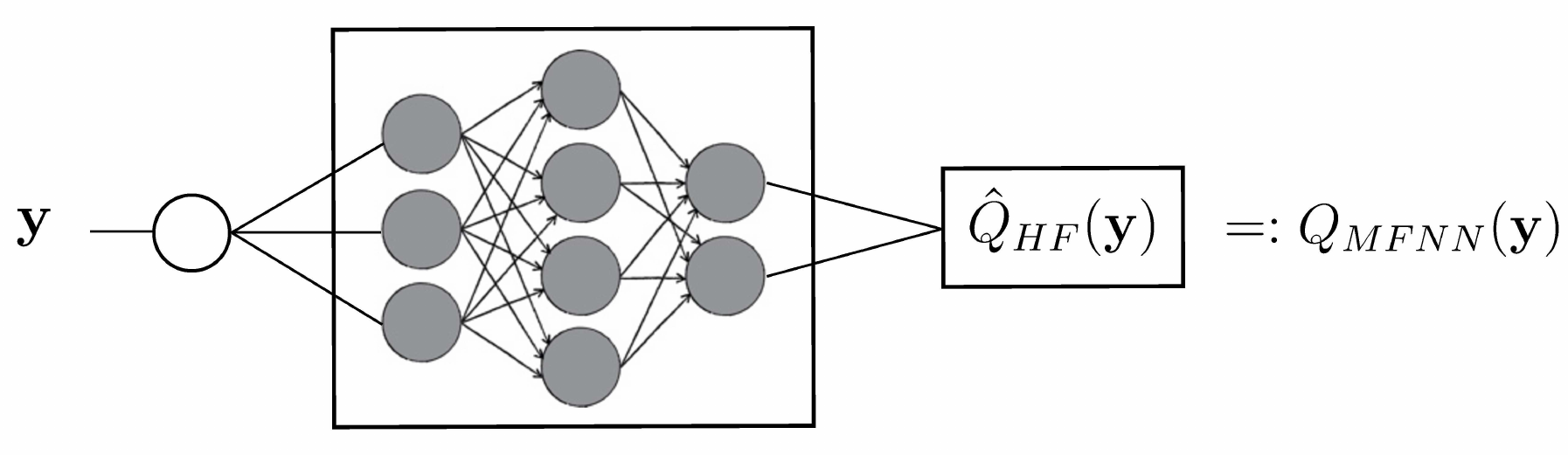}     
\caption{A schematic representation of NN2. It takes the input-output data
  $\{ ( {\bf y}^{(i)}, \hat{Q}_{HF}^{(i)} )\}_{i=1}^{M}$ generated by
  NN1 to learn the high-fidelity quantity $Q_{HF}({\bf y})$.}
\label{NN2_fig}
\end{center}
\end{figure}

The algorithm consists of the following steps. 
\begin{enumerate}
\setlength{\itemsep}{1.5 pt}
\item Generate a set of $M = M_1 + M_2$ realizations of ${\bf y} \in \Gamma$ collected in two disjoint sets: 
$$
Y_I := \{ {\bf y}^{(1)}, \dotsc, {\bf y}^{(M_1)} \} \subset \Gamma, \qquad Y_{II} := \{
  {\bf y}^{(M_1+1)}, \dotsc, {\bf y}^{(M_1+M_2)} \} \subset \Gamma, \qquad Y_I
  \cap Y_{II} = \emptyset.
$$


\item For each ${\bf y}^{(i)} \in Y_I \cup Y_{II}$, with $i=1,
  \dotsc, M$, compute the low-fidelity realizations 
$$
Q_{LF}^{(i)} := M_{DE}^{h_{LF}}
  (\Theta({\bf y}^{(i)})), \qquad i=1,
  \dotsc, M.
$$

\item For each ${\bf y}^{(i)} \in Y_{I}$, with $i=1,
  \dotsc, M_1$, compute the high-fidelity realizations 
\begin{equation}\label{QHF_realizations}
Q_{HF}^{(i)} := M_{DE}^{h_{HF}}
  (\Theta({\bf y}^{(i)})), \qquad i=1,
  \dotsc, M_1.
\end{equation}

\item Using the set of $M_1$ data $\{ ( {\bf y}^{(i)}, Q_{LF}^{(i)},
  Q_{HF}^{(i)}) \}_{i=1}^{M_1}$, construct a neural network NN1 as a surrogate for the correlation
  function \eqref{corr_function}, denoted by $F_{NN}({\bf y},
  Q_{LF}({\bf y}))$; see Figure \ref{NN1_fig}.

\item Use the surrogate model $F_{NN}$ built by NN1 in step 4 and the
  low-fidelity realizations obtained in step 2 to approximate the remaining $M_2$
  high-fidelity quantities. This gives a set of $M$ data  
\begin{equation}\label{QhatHF_realizations}
\hat{Q}_{HF}^{(i)} := \hat{Q}_{HF}({\bf y}^{(i)}) =
\left\{ \begin{array}{l l}
 Q_{HF}({\bf y}^{(i)}), & \quad {\bf y}^{(i)} \in Y_{I} \\
 F_{NN}({\bf y}^{(i)}, Q_{LF}^{(i)}), &  \quad {\bf y}^{(i)} \in Y_{II}
\end{array} \right.,
\end{equation}
to be used in the next step.

\item Using all $M$ high-fidelity data $\{ ( {\bf y}^{(i)},
  \hat{Q}_{HF}^{(i)}) \}_{i=1}^{M}$, construct a neural network NN2 as a surrogate for the high-fidelity quantity $Q_{HF}({\bf
    y})$, denoted by $Q_{MFNN}({\bf y})$; see Figure \ref{NN2_fig}.

\item Generate $N \gg M$ samples of ${\bf y}$ according to the joint PDF, collected in a set of realizations
$$ 
Y := \{ {\bf y}^{(1)}, \dotsc, {\bf y}^{(N)} \}.
$$


\item Approximate the expectation of $Q$ by the sample mean of its
  realizations computed by the multi-fidelity surrogate model NN2 constructed in step 6:
\begin{equation}\label{A_MFNNMC}
{\mathbb E}[Q] \approx {\mathcal A}_{MFNNMC} := \frac1{N}
  \sum_{i=1}^{N} Q_{MFNN} ({\bf y}^{(i)}).
\end{equation}
\end{enumerate}



We note that higher statistical moments may also be computed by taking the sample average of higher powers of the
approximated quantity in the last step of the algorithm.

\begin{remark}\label{remark_MC} 
If we simply take $M=0$, skip steps 1 to 6 above, and compute $N$ 
realizations of $Q$ directly by the high-fidelity model, we
arrive at the classical MC sampling estimation
\begin{equation}\label{A_HFMC}
{\mathbb E}[Q] \approx {\mathcal A}_{HFMC}  := 
\frac1{N}
  \sum_{i=1}^{N} Q_{HF} ({\bf y}^{(i)}).
\end{equation}
\end{remark}

\begin{remark}\label{remark_M1M2}
The two networks (NN1 and NN2) are trained based on two training data sets of
different size: NN1 uses $M_1$ data points, while NN2 uses $M > M_1$
data points. 
This may amount to two networks with different architectures. For instance,
we may consider a simpler architecture (e.g. with a smaller number of
layers/nurons) for NN1 compared to NN2, as the former uses less training data points than the latter. 
While this is true, we do not expect a large difference between the
architectures of NN1 and NN2, as both networks learn the very same
quantities $F$ and ${Q}_{HF}$ and hence need to have comparable
complexities. 
This in turn imposes a restriction on the ratio $r = M_1/M
<1$, that is, $r$ cannot be too small. We need $M_1$ to be large enough to ensure the accuracy of NN1. For instance, the numerical
examples in Section \ref{sec:numerics} use $r=0.25$. 
It is to be noted that in cases when $r \ll 1$, for
instance when we can afford only a few
high-fidelity solves, the
proposed construction may fail to deliver accurate predictions. We
are currently working on a modified construction that enables
accurate predictions even in the case $r \ll 1$, to be presented
elsewhere.
\end{remark}

\begin{remark}\label{Composite_MFNN}
A recently proposed strategy \cite{MFNN:20} also employs NNs to construct multi-fidelity
regression models. The construction in \cite{MFNN:20} has two major
components. First, it uses a combination of a linear and nonlinear
relation between the low- and high-fidelity models: $Q_{HF}({\bf y}) =
F_L ({\bf y}, Q_{LF}({\bf y})) + F_{NL} ({\bf y}, Q_{LF}({\bf
  y}))$. Next, it trains three NNs: one low-fidelity NN for approximating the
  low-fidelity model, and two high-fidelity NNs for learning the linear and
  nonlinear functions $F_{L}$ and $F_{NL}$. Consequently, each
  evaluation of $Q_{HF}({\bf y})$ would need all three networks to be
  evaluated. 
Our construction differs from \cite{MFNN:20} in both
components. 
Instead of the combination of a linear and a nonlinear
function, we employ a general (possibly nonlinear) function $F$ in
\eqref{corr_function}. Such a general function can model both linear
and nonlinear relations between the two models. Consequently, instead
of training and evaluating two networks in \cite{MFNN:20}, one for $F_L$
and one for $F_{NL}$, here we need only one network to learn $F$. 
This choice is motivated by the following observation: a neural
network that captures the (complex) nonlinear features of a function
can also capture its (simple) linear features without the need to add
many extra layers/neurons. 
Moreover, our construction changes the training order proposed in \cite{MFNN:20}
in the following way. While in \cite{MFNN:20} the first network uses
only low-fidelity data and the second and third networks use both
low-fidelity and high-fidelity data, here the first network uses both
low-fidelity and high-fidelity data, and the second network uses only high-fidelity data. This change
in the order of training data makes the ultimate evaluations ``independent'' of the first
network. The first network NN1 is only used to generate additional
training data for the second network NN2, and hence it will be evaluated only a few
times. The two trained networks can therefore be separated, and only the
second network NN2 needs to be evaluated for each evaluation of
$Q_{HF}({\bf y})$. 
As a result, compared to \cite{MFNN:20}, the new construction would save one
network training and two network evaluations. Hence, the new construction may be
superior in terms of both training costs and evaluation costs without sacrificing accuracy. 
\end{remark}

\subsection{Approximation error}
\label{sec:error}

The approximation \eqref{A_MFNNMC} involves two separate estimations: 1) the
estimation of $Q$ by the surrogate model $Q_{MFNN}$ at a set of $N$
realizations of ${\bf y}$, and 2) the estimation of the
integral by a sum of $N$ terms. 
Correspondingly, we may split the total error into two parts:
\begin{equation*}
\varepsilon := | {\mathbb E}[Q({\bf y})] - {\mathcal A}_{MFNNMC}|  
\leq 
\underbrace{| {\mathbb E}[Q({\bf y})] - {\mathbb E}[Q_{MFNN}({\bf y})] |}_{\varepsilon_{I}} 
+ 
\underbrace{| {\mathbb E}[Q_{MFNN}({\bf y})] -  {\mathcal A}_{MFNNMC}|}_{\varepsilon_{II}}.
\end{equation*}

The first error term $\varepsilon_{I}$ is deterministic and
corresponds to the error in the approximation of $Q({\bf y})$ by the
multi-fidelity neural network surrogate model $Q_{MFNN}({\bf y})$. 
Two sources contribute to this error: 
1) the structure of NN2, including the number of layers and neurons,
activation functions, and the regularization and
optimization techniques and parameters used to train the network, and 
2) the quantity and quality of the data used to train NN2, i.e. the
number $M$ and location of the input points $\{ {\bf y}^{(i)}
\}_{i=1}^M \in Y_{I} \cup Y_{II} \subset \Gamma$ and the accuracy of the approximate realizations $\{
\hat{Q}_{HF}^{(i)} \}_{i=1}^{M}$. 
The accuracy of the output data $\{ \hat{Q}_{HF}^{(i)} \}_{i=1}^{M}$ given in
\eqref{QhatHF_realizations} is determined by the error in approximating
the QoI by the high-fidelity model, i.e. the
error in $Q({\bf y}) \approx Q_{HF}({\bf y})$ that satisfies
\eqref{HF_error}, and the error in approximating the high-fidelity
quantity by NN1, i.e. the error in $Q_{HF}({\bf y}) \approx F_{NN}({\bf
y}, Q_{LF}({\bf y}))$. The latter in turn depends on two factors: i) the structure
of NN1, and ii) the quantity and quality of the data used to train NN1, i.e. the
number $M_1$ and location of the input points $\{ {\bf y}^{(i)}
\}_{i=1}^{M_1} \in Y_{I} \subset \Gamma$ and the accuracy of the
approximate realizations $\{ Q_{HF}^{(i)} \}_{i=1}^{M_1}$, which again
satisfies \eqref{HF_error}. 
In summary, the first error term $\varepsilon_{I}$ depends on the
structures of the two trained neural networks NN1 and NN2, the choice of the sets $Y_{I}$ and $Y_{II}$, and the error in the approximation of $Q$ by
the high-fidelity model satisfying \eqref{HF_error}. 
As mentioned in Section \ref{sec:NN}, the precise
dependence of the output error of a trained network on the network
structure and the choice of input training data is still an
open problem. However, thanks to the universal approximation theorem (see e.g. \cite{Hornik_etal:89}), we may assume that there exist neural networks NN1 and NN2
that deliver predictions $Q_{MFNN}$ that are as accurate as the output training data $Q_{HF}$. More precisely, we assume that with a
proper selection of network structure and input training data, we
have
\begin{equation}\label{eps1}
\varepsilon_{I} \le C \, h_{HF}^q.
\end{equation}
It is to be noted that while the universal
approximation theorem states that there exists a network that can
approximate a function within a desired tolerance, it does not tell
how to construct such a network. Hence, we cannot a priori guarantee that the
networks that we train would produce an error satisfying
\eqref{eps1}. 
In practice, we
may need more validation data sets to reduce the possibility that
\eqref{eps1} does not hold.

The second error term $\varepsilon_{II}$ is a quadrature error due
to approximating the integral by a sum. In the particular case of MC
sampling, since ${\mathcal A}_{MFNNMC}$ is a statistical term,
$\varepsilon_{II}$ is referred to as the statistical error. 
By the central limit
theorem (see e.g. \cite{MC}), we know that for large $N$ the
distribution of the term $ {\mathbb E}[Q_{MFNN}({\bf y})] -  {\mathcal A}_{MFNNMC}$ will
approach a normal distribution with mean zero and variance
$\sqrt{{\mathbb V}[Q_{MFNN}] / N}$. 
Consequently, the error $\varepsilon_{II}$ satisfies
\begin{equation}\label{eps2}
P \left( \varepsilon_{II} \le c_{\alpha} \, \sqrt{\frac{{\mathbb
      V}[Q_{MFNN}]}{N}} \, \right) \rightarrow 2 \phi (c_{\alpha}) - 1, \qquad
\text{as} \ \ N \rightarrow \infty,
\end{equation}
where $P$ is a probability measure, and $\phi (c_{\alpha}) =
\frac{1}{\sqrt{2 \pi}}
\int_{-
  \infty}^{c_{\alpha}} \exp(- \tau^2/2) d\tau$ is the cumulative
density function of a standard normal random variable evaluated at a
given confidence level $c_{\alpha}>0$. 
Clearly, the larger the confidence level, the higher the probability that
the inequality $\varepsilon_{II} \le c_{\alpha} \, \sqrt{{\mathbb
    V}[Q_{MFNN}] / N}$ holds. 

We often need the approximation ${\mathcal A}_{MFNNMC}$ to be accurate within a
desired small tolerance $\varepsilon_{\footnotesize{\text{TOL}}} > 0$ and with a
pre-assigned small failure probability $\alpha \in (0,1)$, that is,
$$
P( \varepsilon  \le \varepsilon_{\footnotesize{\text{TOL}}} ) = 1- \alpha.
$$ 
To achieve this, we may first split the tolerance between the deterministic
and statistical errors by introducing a splitting parameter $\theta
\in (0,1)$ and require that
\begin{equation}\label{two-inequalities}
\varepsilon_{I} \le (1- \theta)
\varepsilon_{\footnotesize{\text{TOL}}}, \qquad P( \varepsilon_{II}
\le \theta \, \varepsilon_{\footnotesize{\text{TOL}}} ) = 1- \alpha. 
\end{equation}
From \eqref{eps1} and the first inequality in \eqref{two-inequalities} we obtain
$h_{HF}$ by requiring $C h_{HF}^q \le (1- \theta) \,
\varepsilon_{\footnotesize{\text{TOL}}}$. Moreover, comparing \eqref{eps2}
and the second
inequality in \eqref{two-inequalities}, we set $2 \phi (c_{\alpha}) -
1 = 1-\alpha$ and obtain the confidence level $c_{\alpha} = \phi^{-1}(1-\alpha/2)$. Then, the number of
samples $N$ will be obtained by requiring 
$c_{\alpha} \, \sqrt{{\mathbb
    V}[Q_{MFNN}] / N} \le \theta \,
\varepsilon_{\footnotesize{\text{TOL}}}$. Finally, we get
$$
h_{HF} \propto \varepsilon_{\footnotesize{\text{TOL}}}^{1/q}, \qquad N
\propto \varepsilon_{\footnotesize{\text{TOL}}}^{-2}.
$$

It is to be noted that \eqref{two-inequalities} guarantees that
$P(\varepsilon  > \varepsilon_{\footnotesize{\text{TOL}}} ) \le
\alpha$, which is more conservative than what we need: $P(\varepsilon
> \varepsilon_{\footnotesize{\text{TOL}}} ) = \alpha$. In order to
avoid too conservative considerations, we should select $h_{HF}$ so that the deterministic error
$\varepsilon_I$ is as close as possible to $ (1 - \theta)
\varepsilon_{TOL}$.

\subsection{Computational cost}
\label{sec:cost}

Let $W_{LF}$ and $W_{HF}$ denote the computational cost of evaluating
$Q_{LF}$ and $Q_{HF}$ in \eqref{bi-fidelity} at a single realization of
${\bf y}$, respectively. 
Let further $W_{Ti}$ and $W_{Pi}$, with $i=1,2$, be the training cost
and the cost of evaluating NN$i$ at a single realization of
${\bf y}$, respectively. The total cost of computing the estimator
${\mathcal A}_{MFNNMC}$ in \eqref{A_MFNNMC} will then be
\begin{equation}\label{MFNNMC_cost}
W_{MFNNMC} = M \, W_{LF} + M_1 \, W_{HF} + W_{T1} + M_2 \, W_{P1} +W_{T2} + N \, W_{P2}.
\end{equation}
%
We will now discuss the efficiency of the proposed sampling method and
compare its cost with the cost of the classical high-fidelity
MC estimator ${\mathcal A}_{HFMC}$ formulated in \eqref{A_HFMC}, which is
\begin{equation}\label{HFMC_cost}
W_{HFMC} = N \, W_{HF}.
\end{equation}

In many problems, we often need to obtain accurate predictions within small
tolerances $\varepsilon_{\footnotesize{\text{TOL}}} \ll 1$. To achieve
this within the framework of MC sampling, we usually need a very large
number of samples $N \propto \varepsilon_{\footnotesize{\text{TOL}}}^{-2} \gg 1$. This is because the statistical error in
MC sampling is proportional to $N^{-1/2}$ and hence decays very slowly
as $N$ increases, as discussed above. The classical high-fidelity MC sampling would
therefore be very expensive as it requires a very large number $N$ of
expensive high-fidelity problems to be solved. For instance, suppose $W_{HF} \propto h_{HF}^{- \gamma}$, where $\gamma>0$ is related
to the time-space dimension of the ODE/PDE problem and the discretization
technique used to solve the problem. Then, noting $h_{HF} \propto
\varepsilon_{\footnotesize{\text{TOL}}}^{1/q}$, the cost of MC
sampling \eqref{HFMC_cost} reads
\begin{equation}\label{HFMC_cost_tol}
W_{HFMC} \propto \varepsilon_{\footnotesize{\text{TOL}}}^{-(2+ \gamma/q)}.
\end{equation}
In the new sampling
algorithm proposed here, however, we may achieve dramatic savings in computational
cost by replacing many (i.e. $N$) expensive high-fidelity ODE/PDE solves needed in MC sampling
with the same number of fast neural network computations and only a few (i.e. $M_1$) high-fidelity solves. 
%
%
For this to work, 
we would need: 1) the number $M$ of required training data to be much
less than the number $N$ of required MC samples; 2) the cost $W_{P2}$ of
evaluating NN2 to be much less than the cost $W_{HF}$ of a complex ODE/PDE
solve; and 3) the total training cost $W_{T1}+W_{T2}$ to be negligible
compared to the cost of solving $N$ complex ODE/PDE problems. 
A few observations addressing each of the above three requirements
follow. See also Section \ref{sec:numerics} for numerical
verifications. 
\begin{itemize}[topsep=6pt,leftmargin=3.85\labelsep]
\setlength\itemsep{-.05em}
\item While there is currently no results
on the dependency of the number of required training data on
accuracy, we have observed through numerical experiments that
$M$ depends mildly on the tolerance, i.e. $M \propto
\varepsilon_{\footnotesize{\text{TOL}}}^{-p}$, where $p \in [0,1)$ is
small. A similar observation has also been made in
\cite{Afroja:19}. A comparison between $M \propto
\varepsilon_{\footnotesize{\text{TOL}}}^{-p}$ and $N \propto
\varepsilon_{\footnotesize{\text{TOL}}}^{-2}$ 
suggests that the smaller the tolerance, the
larger the ratio $N/M$. Consequently, the term $M \, W_{LF} + M_1 \,
W_{HF}$ in \eqref{MFNNMC_cost} will be proportional to
$\varepsilon_{\footnotesize{\text{TOL}}}^{-(p+\gamma/q)}$, which is
much less than $N \, W_{HF}$ and hence negligible at small tolerance levels.

\item The cost of evaluating a neural network depends on the network
  architecture, i.e. the number of hidden layers and neurons and the
  type of activation functions. The
  evaluation of a network mainly involves simple matrix-vector
  operations and the application of activation functions. Crucially,
  the evaluation cost is independent of the complexity of the ODE/PDE problem. In
  particular, for QoIs that can be well approximated by a network with
  a few hidden layers and hundreds of neurons,
  we would get $W_{Pi} \ll W_{HF}$, with $i=1,2$. Indeed, in such cases, the more complex the high-fidelity
  problem and the more expensive computing the high-fidelity quantity, the larger the ratio
  $W_{HF} / W_{Pi}$.

\item The cost of training the two networks depends on the number of
  training data, network structure, and the cost of solving their corresponding
  optimization problems. Importantly, it is a one-time cost. In particular, we have
  observed through numerical experiments that as the tolerance
  decreases and $N$ increases, this cost becomes negligible.
\end{itemize}
In summary, for the ODE/PDE problems and QoIs where the above three
requirements hold, we may obtain for small tolerances:
\begin{equation}\label{MFNNMC_cost_tol}
W_{MFNNMC} \propto \varepsilon_{\footnotesize{\text{TOL}}}^{-\max(2, \, p+\gamma/q)}.
\end{equation} 
Comparing \eqref{HFMC_cost_tol} and \eqref{MFNNMC_cost_tol}, we
consider two cases: 1) if $p + \gamma/q \le 2$, then $W_{MFNNMC} \propto
\varepsilon_{\footnotesize{\text{TOL}}}^{-2}$ will be much smaller
than $W_{HFMC}$, and 2) if $p + \gamma/q > 2$, then $W_{MFNNMC} \propto
\varepsilon_{\footnotesize{\text{TOL}}}^{-(p + \gamma/q)}$ will be much smaller
than $W_{HFMC}$ as long as $p<2$. 
Importantly, the smaller the tolerance $\varepsilon_{\text{TOL}}$, the more gain in
computational cost over MC sampling.

\section{Numerical Examples}
\label{sec:numerics}

In this section we present two numerical examples: an ODE problem and
a PDE problem. 
All codes are written in Python and run on a single CPU in order to
have a fair comparison between the proposed method and MC sampling. We use Keras
\cite{chollet2015keras}, which is an open-source neural-network
library written in Python, to construct the neural networks. 
It is to be noted that all CPU times are  measured by {\tt time.clock()} in Python. 

\subsection{An ODE problem}

Consider the following parametric initial value problem (IVP)
\begin{equation}\label{IVP}
\begin{array}{ll}
u_t(t,y) + 0.5 u(t,y) = f(t,y), & \qquad t \in [0,T], \quad y \in
  \Gamma \subset {\mathbb R}, \\
u(0,y)=g(y).& 
\end{array}
\end{equation}
where $y \in \Gamma$ is a uniformly distributed random
variable on $\Gamma = [-1,1]$, and the force term $f$ and the initial
data $g$ are so that the exact solution to the IVP \eqref{IVP} is
$$
u(t,y) = 0.5+2 \sin(12 y)+6 \sin(2t) \sin(10y) (1+2 y^2).
$$
Our goal is to approximate the expectation ${\mathbb E}[Q(y)]$,
where $Q(y) =| q(y) | =| u(T,y)|$ with $T=100$, by the multi-fidelity
estimator ${\mathcal A}_{MFNNMC}$ in \eqref{A_MFNNMC} and compare its
performance with the high-fidelity MC estimator ${\mathcal A}_{HFMC}$
in \eqref{A_HFMC}. We use the closed form of solution to measure
errors and will compare the cost of the two methods subject to the
same accuracy constraint.

Suppose that we have the second-order accurate Runge-Kutta (RK2)
time-stepper as the deterministic solver to compute realizations of
$q_{LF}(y)$ and $q_{HF}(y)$ using time steps $h_{LF}$ and $h_{HF}$,
respectively. 
Consider the relative error
in the approximation
$$
\varepsilon_{\text{rel}} := | {\mathbb E}[Q(y)] - {\mathcal
  A}|  / | {\mathbb E}[Q(y)] |,
$$
where the estimator ${\mathcal A}$ is either ${\mathcal A}_{MFNNMC}$
or ${\mathcal A}_{HFMC}$. 
Given a $1\%$ failure probability ($\alpha=
0.01$) and a decreasing sequence of tolerances
$\varepsilon_{\footnotesize{\text{TOL}}} = 10^{-2}, 10^{-3}, 10^{-4}$,
a simple error analysis similar to the analysis in
Section \ref{sec:error} and verified by numerical computations gives
the minimum number of realizations $N$ and the maximum time step
$h_{HF}$ for the high-fidelity model required to achieve $P( \varepsilon_{\text{rel}}  \le
\varepsilon_{\footnotesize{\text{TOL}}} ) = 0.99$. We choose a fixed time step $h_{LF}=0.5$ for the
low-fidelity model at all tolerance levels. Table
\ref{Ex1_table1} summarizes the numerical parameters $(N, h_{HF}, h_{LF})$ and the CPU time
of evaluating single realizations of $q_{LF}$ and $q_{HF}$.   
\begin{table}[!h]
        \centering
        \caption{Required number of realizations and time steps to
        achieve $P( \varepsilon_{\text{rel}}  \le
\varepsilon_{\footnotesize{\text{TOL}}} ) = 0.99$.}
\medskip
        \label{Ex1_table1}
        \begin{tabular}{|c|c|c|c|c|c|}
\hline
            $\varepsilon_{\footnotesize{\text{TOL}}}$ & $N$ & $h_{HF}$ &
          $W_{HF}$ & $h_{LF}$ & $W_{LF}$ \\
            \hline
            \hline
             $10^{-2}$ & $1.35 \times 10^{5}$ & 0.1     & $2.24 \times
                                                          10^{-4}$ & 0.5 & $4.36 \times
                                                          10^{-5}$\\
             $10^{-3}$ & $1.35 \times 10^{7}$ & 0.025 & $7.21 \times
                                                          10^{-4}$ & 0.5 & $4.36 \times
                                                          10^{-5}$\\
             $10^{-4}$ & $1.35 \times 10^{9}$ & 0.01   & $2.20 \times
                                                          10^{-3}$ & 0.5 & $4.36 \times
                                                          10^{-5}$\\
            \hline
        \end{tabular}
\end{table}
%

%

\begin{figure}[!h]
\begin{center}
\includegraphics[width=0.85\linewidth]{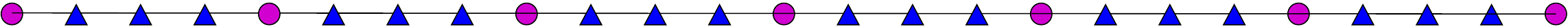}     
\caption{A schematic representation of the selected points. Magenta
  circles represent the points in $Y_{I}$, where only the
  high-fidelity quantity is computed. The blue triangles represent
  the points in $Y_{II}$. The low fidelity quantity is computed at all
points in $ Y_I \cup Y_{II}$.}
\label{NN_points}
\end{center}
\end{figure}

Following the algorithm in Section \ref{sec:algorithm}, we first
generate a set of $M = M_1 + M_2$ points $y^{(i)} \in [-1,1]$, with
$i=1, \dotsc, M$, collected into two disjoint sets $Y_I$ and $Y_{II}$. Figure \ref{NN_points} shows a
schematic representation of the selection of points. 
We choose the points to be
uniformly placed on the interval $[-1,1]$. We then select every
4th point to be in the set $Y_I$ (magenta circles) and collect the rest of the points in
the set $Y_{II}$ (blue triangles). This implies $M_2 \approx 3 \, M_1$, meaning that we
will need to compute the quantity $q(y)$ by the high-fidelity model,
that is RK2 using time step $h_{HF}$, at only a quarter of points
$M/4$. The number of points $M$ will be chosen based on the desired
tolerance, slightly increasing as the tolerance decreases. 


\begin{figure}[!h]
\vspace{-0.2cm}
\begin{center}
\includegraphics[width=0.55\linewidth]{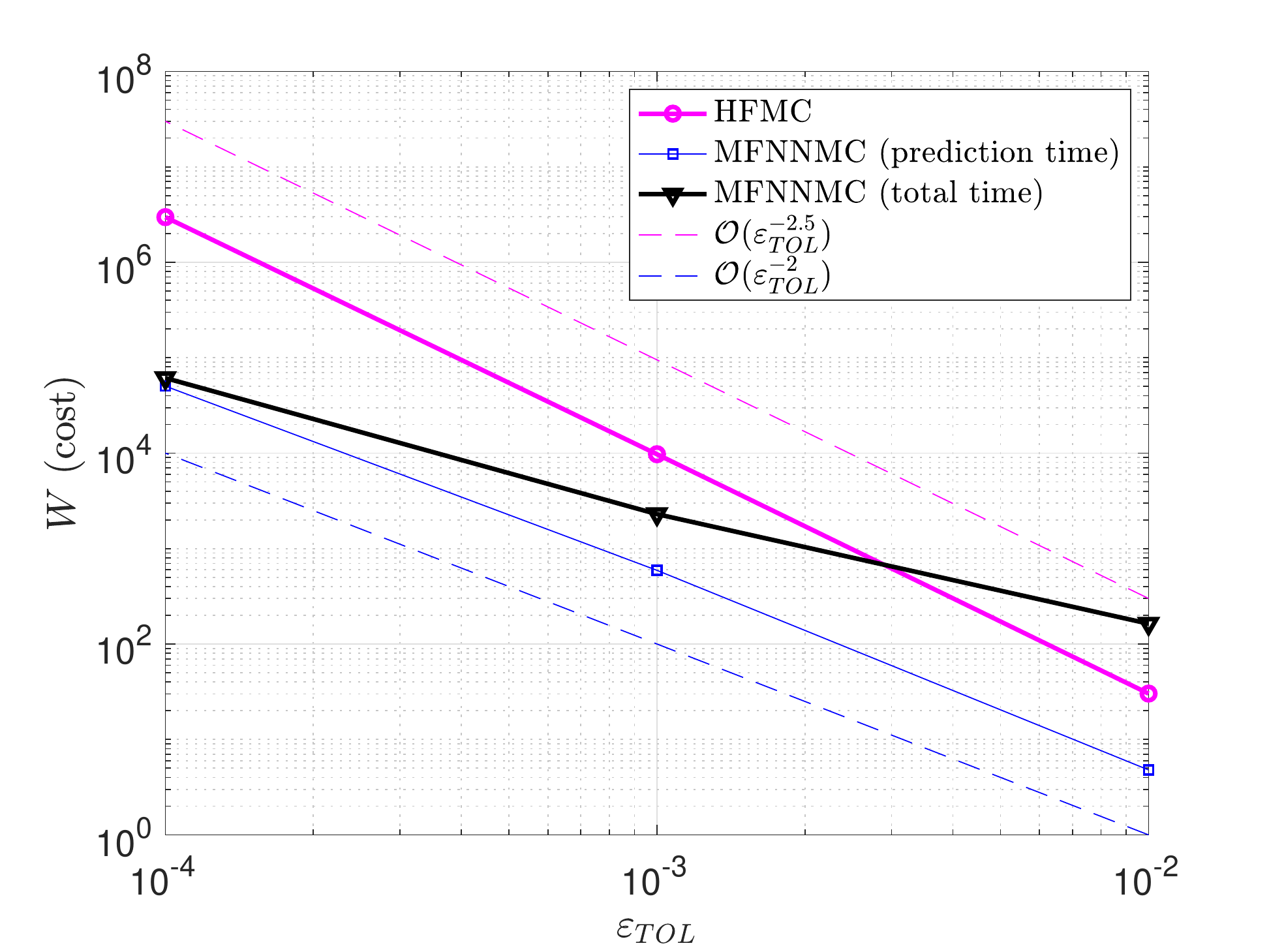}     
\vspace{-0.5cm}
\caption{CPU time versus tolerance. For large tolerances the
  training cost is dominant, making the cost of
  MFNNMC more than the cost of HFMC. However, as tolerance decreases,
  the training cost becomes negligible and the cost of
  MFNNMC approaches ${\mathcal
  O}(\varepsilon_{\footnotesize{\text{TOL}}}^{-2})$.}
\vspace{-.5cm}
\label{Ex1_conv}
\end{center}
\end{figure}

\begin{figure}[!h]
\vspace{-0.25cm}
\begin{center}
\includegraphics[width=0.55\linewidth]{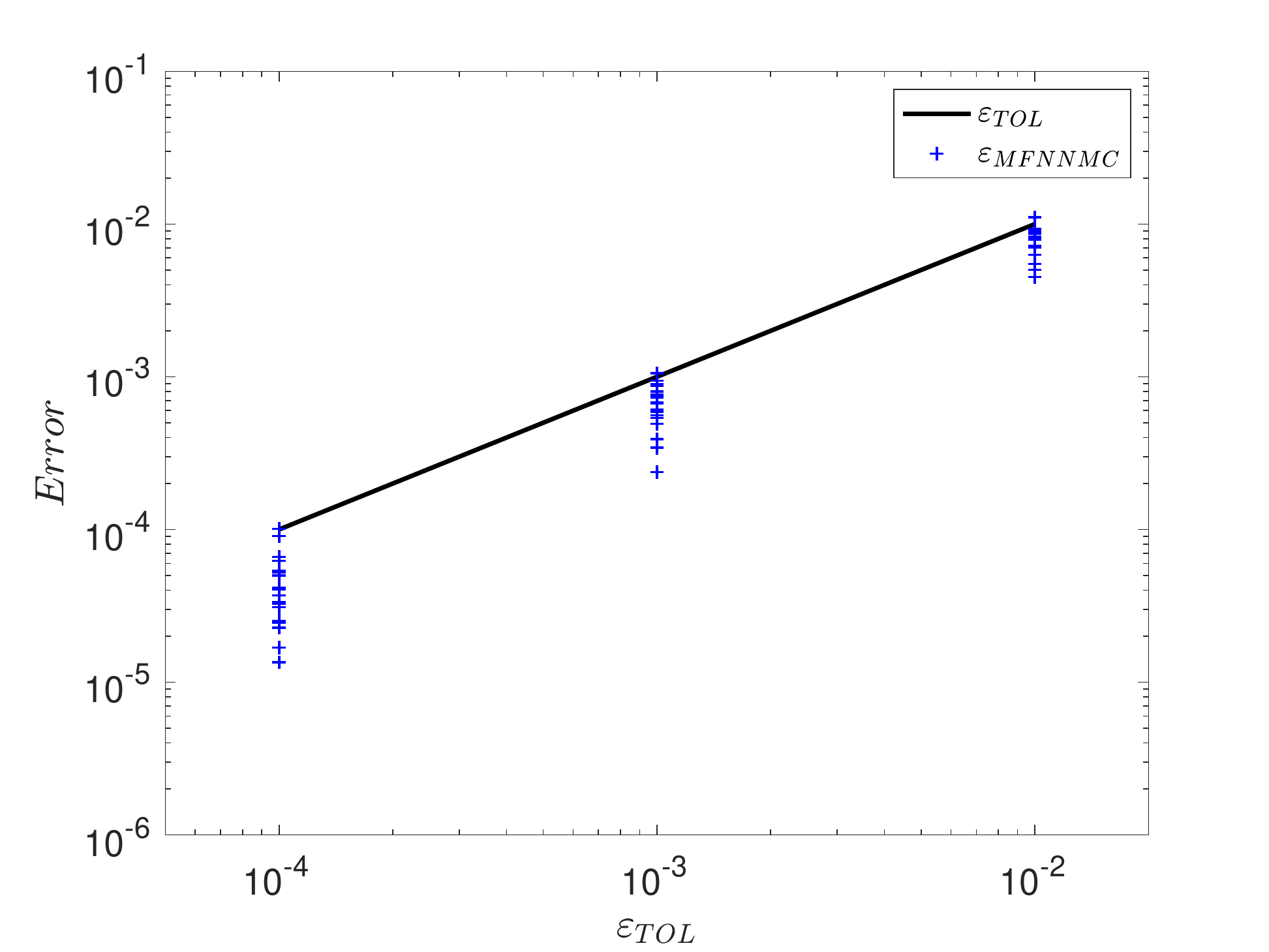}     
\vspace{-0.5cm}
\caption{Relative error as a function of tolerance, verifying that the
  tolerance is met with 1$\%$ failure probability. The ``+'' markers
  correspond to 20 simulations at each tolerance level.}
\label{Ex1_tols}
\end{center}
\end{figure}

We will use the same architecture for the two networks NN1 and NN2 and
keep them fixed at all tolerances. Precisely, we choose feed-forward networks with
4 hidden layers, where each layer contains 20 neurons. We use ReLU
activation function for the hidden layers and the identity activation
function for the output layer of both networks. It is to be noted that
NN1 has two input neurons, while NN2 has one input neuron. Both networks
have one output neuron. For the training process, we employ the
quadratic cost function (or the mean squared error) and use the Adam
optimization technique with a fixed learning rate $\eta = 0.002$. 
We note that instead of a fixed learning rate, one may alternatively
split the available $M$ data points into training and validation sets
and tune the learning rate parameter. For this ODE example, we do not
use any validation set as the choice $\eta = 0.002$ produces
satisfactory results both in terms of efficiency and accuracy. 
We also do not use any regularization technique. 
Table \ref{Ex1_table2}
summarizes the number of training data $M=M_1 + M_2$, the number of epochs $N_{\text{epoch}}$, batch size
$N_{\text{batch}}$, and the CPU time of training and evaluating the two
networks for different tolerances. 
We observe that the number of training data satisfies $M
\propto \varepsilon_{\footnotesize{\text{TOL}}}^{-p}$ with $p=0.5$. 
It is also to be noted that while using a different architecture for each
network and other choices of network parameters
(e.g. number of layers/ neurons and learning rates) may give
more efficient networks, the selected architectures and parameters here,
following the general guidelines in \cite{Bengio:12,GoodBengCour16}, produce
satisfactory results in terms of efficiency (see Figure \ref{Ex1_conv}) and
accuracy (see Figure \ref{Ex1_tols}). 
\begin{table}[!h]
       \centering
        \caption{The number of training data and training and
          evaluation time of the two networks.}
\medskip
        \label{Ex1_table2}
\resizebox{\textwidth}{!}{
\begin{tabular}{|c||c|c||c|c|c|c||c|c|c|c|}
\hline
& & &\multicolumn{4}{c||}{NN1}&\multicolumn{4}{c|}{NN2}\\
\cline{4-11}
$\varepsilon_{\footnotesize{\text{TOL}}}$ & $M_1$ & $M_2$ &
                                                            $N_{\text{epoch}}$
                                          & $N_{\text{batch}}$ &
                                                                 $W_{T_1}$
                                          & $W_{P_1}$ &
                                                            $N_{\text{epoch}}$
                                          & $N_{\text{batch}}$ &
                                                                 $W_{T_2}$
                                          & $W_{P_1}$\\
\hline\hline
$10^{-2}$ & 61   & 180   &100& 10& 10.98 & $3.57 \times
                                                          10^{-5}$ &
                                                                     1800&
                                                                           40&
                                                                               147.44 & $3.55 \times 10^{-5}$   \\
$10^{-3}$ & 201 & 600   &1000& 30& 81.13 & $7.77 \times
                                                          10^{-5}$ & 5000& 80& 1623.41 & $4.38 \times 10^{-5}$   \\
$10^{-4}$ & 801 & 2400 &5000& 100& 768.30 & $8.39 \times
                                                          10^{-5}$ &30000& 100& 10096.63 & $3.73 \times 10^{-5}$  \\
\hline
\end{tabular}}
\end{table}

\begin{figure}[!h]
\begin{center}
\includegraphics[width=0.55\linewidth]{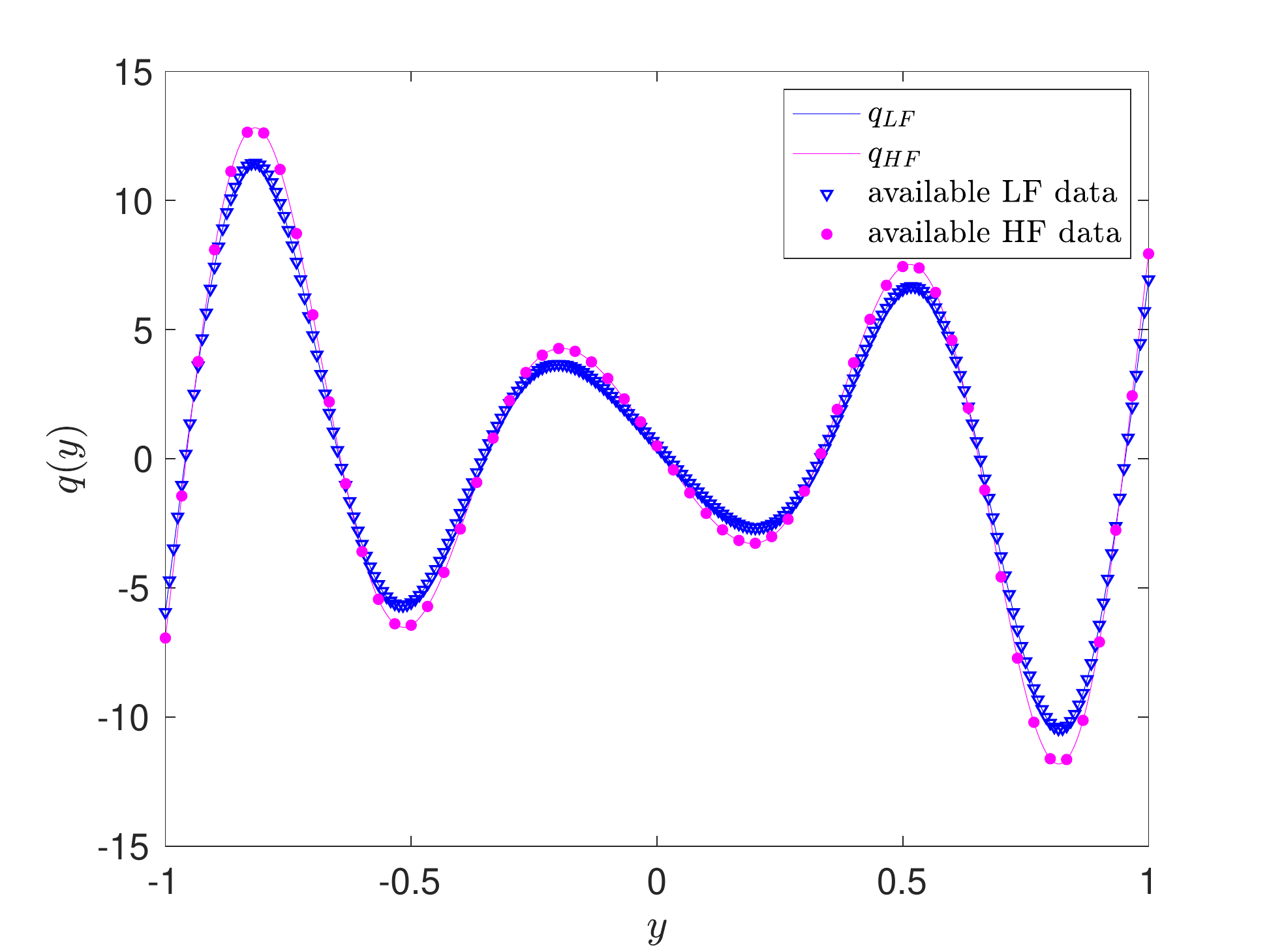}     
\caption{The low-fidelity and high-fidelity quantities versus
  $y\in[-1,1]$ (solid lines) and the
  available data (markers) in the case $\varepsilon_{\footnotesize{\text{TOL}}} =
  10^{-2}$. There are $M_1 = 61$ high-fidelity and $M=241$
  low-fidelity data points, represented by circles and triangles, respectively.}
\label{Ex1_tol1_points}
\end{center}
\end{figure}

\begin{figure}[!h]
\vspace{-.1cm}
\center
\subfigure{\includegraphics[width=8.1cm]{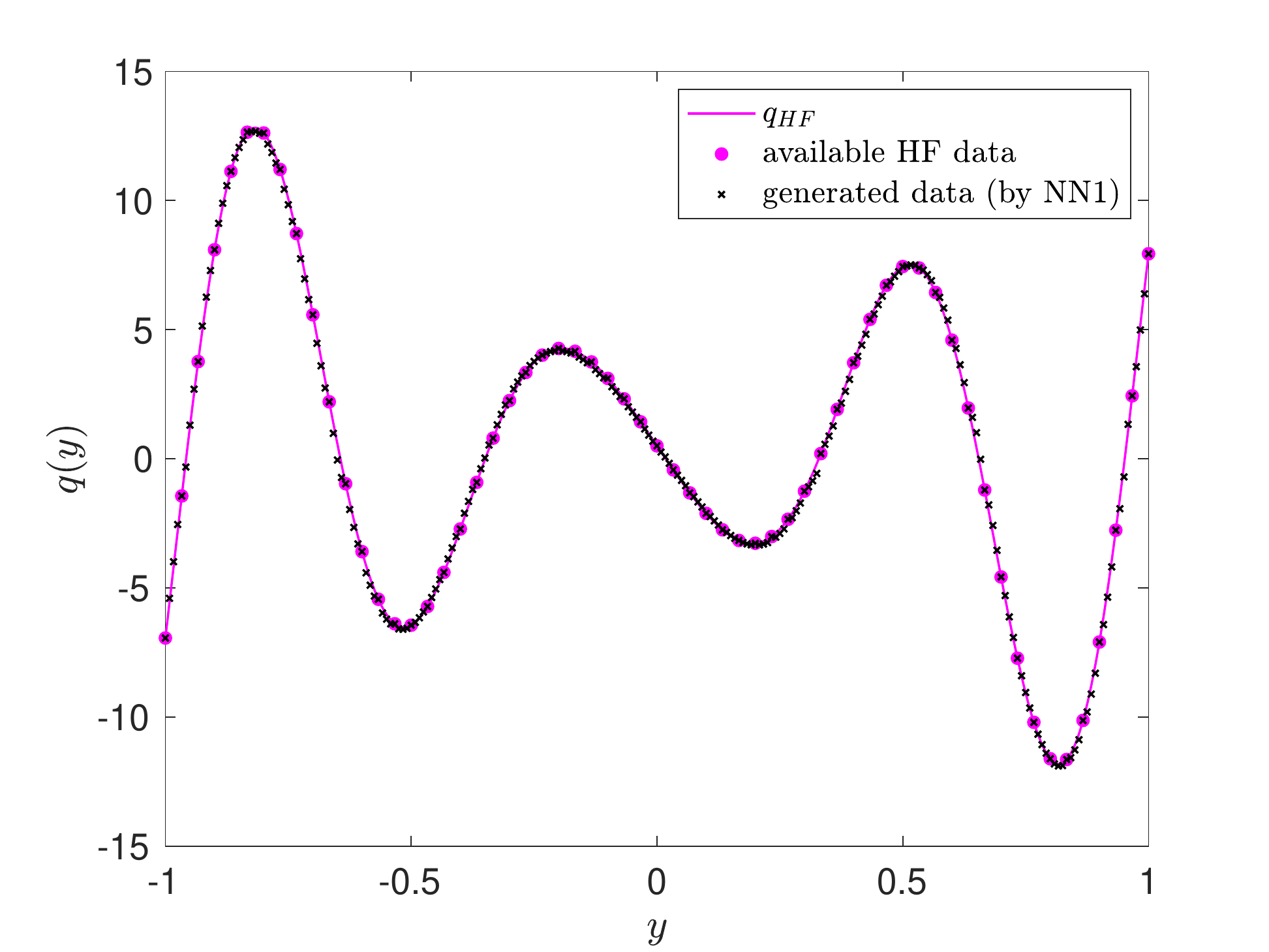}}
 \hskip .1cm
\subfigure{\includegraphics[width=8.1cm]{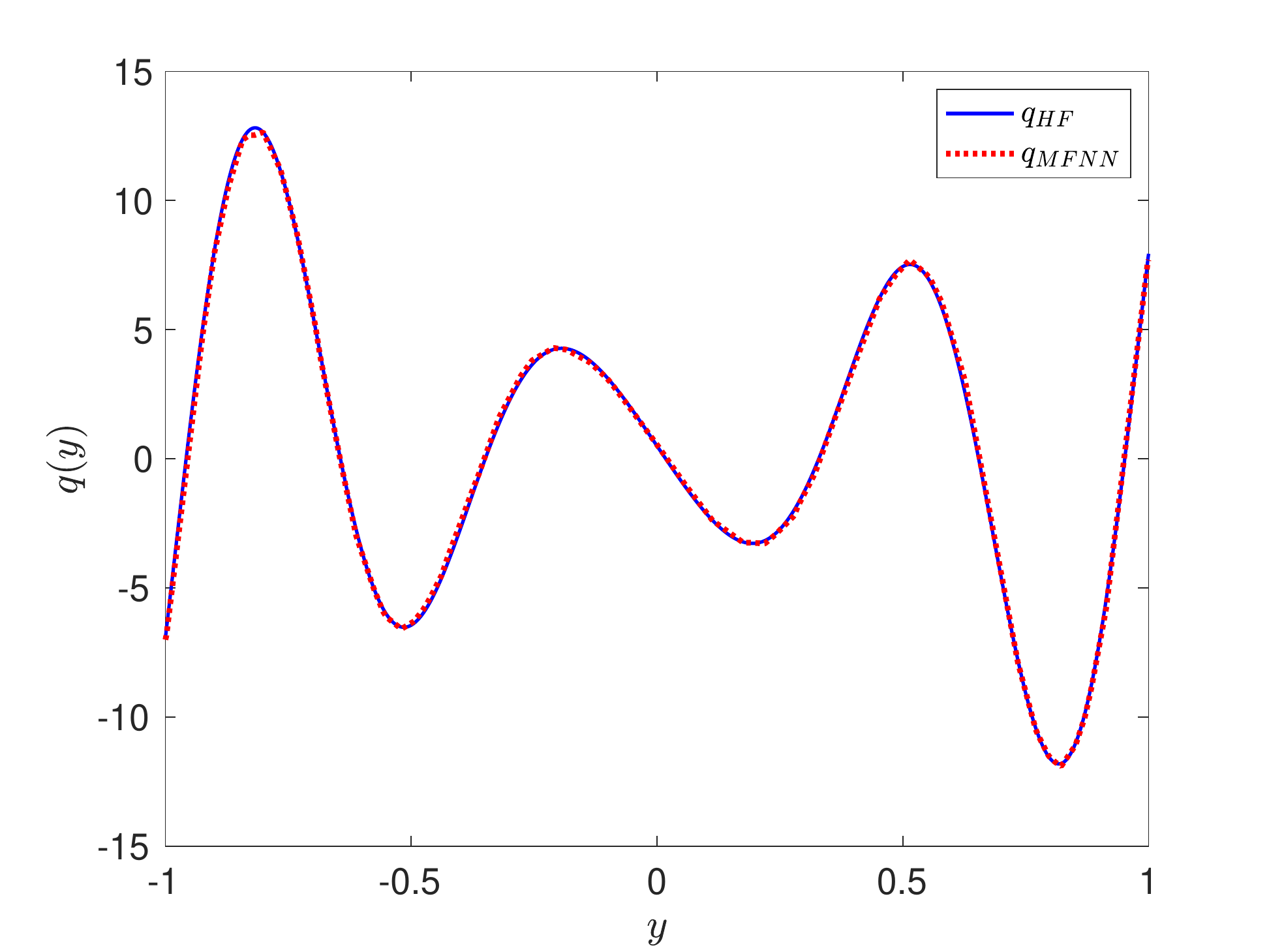}}
\vspace{-.1cm}
    \caption{Outputs of the trained networks for $\varepsilon_{\footnotesize{\text{TOL}}} =
  10^{-2}$. Left: generated data by the trained network NN1. Right:
      predicted quantity by the trained network NN2.}
    \label{Ex1_toal1_predictions}
 \end{figure}

Figure \ref{Ex1_tol1_points} shows the low-fidelity and high-fidelity
quantities versus $y \in [-1,1]$ (solid lines) and the data
(circle and triangle markers) available in the case $\varepsilon_{\footnotesize{\text{TOL}}}=10^{-2}$. 
Figure \ref{Ex1_toal1_predictions} (left) shows the generated
high-fidelity data by the trained network NN1, and Figure
\ref{Ex1_toal1_predictions} (right) shows the predicted high-fidelity quantity by the
trained network NN2 for tolerance $\varepsilon_{\footnotesize{\text{TOL}}}=10^{-2}$. 

Figure \ref{Ex1_conv} shows the CPU time as a function of
tolerance. The computational cost of classical MC sampling is
${\mathcal O}(\varepsilon_{\footnotesize{\text{TOL}}}^{-2.5})$, following \eqref{HFMC_cost_tol} and noting that the order of
accuracy of RK2 is $q=2$ and the time-space dimension of the problem
is $\gamma = 1$. 
On the other hand, if we only consider the prediction time of the
proposed multi-fidelity method, excluding the training costs, the cost of the proposed method is proportional to $\varepsilon_{\footnotesize{\text{TOL}}}^{-2}$ which
is much less than the cost of MC sampling. When adding the training
costs, we observe that although for large tolerances the training cost
is large, as the tolerance decreases the training costs become
negligible compared to the total CPU time. Overall, the cost of the
proposed method approaches ${\mathcal
  O}(\varepsilon_{\footnotesize{\text{TOL}}}^{-2})$ as tolerance
decreases, and hence, the smaller the tolerance, the more gain in
computational cost when employing the proposed method over MC
sampling. This can also be seen by \eqref{MFNNMC_cost_tol} where
$\max(2, p + \gamma/q) = \max(2, 0.5 + 0.5) = 2$. 

Finally, Figure \ref{Ex1_tols} shows the relative error as a function
of tolerance for the proposed method, verifying that the tolerance
is met with 1$\%$ failure probability.

\subsection{A PDE problem}

Consider the following parametric initial-boundary value problem
(IBVP)
\begin{equation}\label{IBVP}
\begin{array}{ll}
u_{tt}(t,{\bf x},{\bf y})- \Delta_{\bf x} u(t,{\bf x},{\bf y}) = f(t,{\bf x},{\bf y}),  & \ \ \
                                                                    (t, {\bf x},{\bf y}) \in [0,T] \times D \times \Gamma,\\
u(0, {\bf x},{\bf y}) = g_1({\bf x},{\bf y}), \ \ u_t(0, {\bf x},{\bf y}) = g_2({\bf x},{\bf y}), &  \ \ \ (t,{\bf x},{\bf y})
                                        \in \{ 0 \} \times D \times \Gamma, \\
u(t, {\bf x},{\bf y}) = g_b(t,{\bf x},{\bf y}), &  \ \ \ (t,{\bf x},{\bf y})
                                        \in [0,T] \times \partial D \times \Gamma,
\end{array}    
\end{equation}
where $t \in [0,T]$ is the time, ${\bf x} = (x_1, x_2) \in D$ is the vector of
spatial variables on a square domain $D=[-1,1]^2$, and ${\bf y} = (y_1, y_2) \in \Gamma$ is a vector of two uniformly distributed random
variables on $\Gamma = [10,11]\times[4,6]$. We select the force term $f$ and the initial-boundary
data $g_1, g_2, g_b$ so that the exact solution to the IBVP \eqref{IBVP} is
$$
u(t,{\bf x},{\bf y}) = \sin(y_1 \, t - y_2 \, x_1) \, \sin(y_2 \, x_2).
$$
Our goal is to approximate the expectation ${\mathbb E}[Q(y)]$,
where $Q({\bf y}) =| q({\bf y}) | =| u(T,{\bf x}_Q, {\bf y})|$ with
$T=30$ and ${\bf x}_Q = (0.5,0.5)$, by the multi-fidelity
estimator ${\mathcal A}_{MFNNMC}$ in \eqref{A_MFNNMC} and compare its
performance with the high-fidelity MC estimator ${\mathcal A}_{HFMC}$
in \eqref{A_HFMC}. We use the closed form of solution to measure
errors and will compare the cost of the two methods subject to the
same accuracy constraint.

Suppose that we have a second-order accurate (in both time and space) finite difference scheme as the deterministic solver to compute realizations of
$q_{LF}({\bf y})$ and $q_{HF}({\bf y})$ using a uniform grid with grid
lengths $h_{LF}$ and $h_{HF}$,
respectively. We use the time step $\Delta t = h/2$ to ensure
stability of the numerical scheme, where the grid length $h$
is either $h_{LF}$ or $h_{HF}$, depending on the level of fidelity. 
Consider the absolute error
in the approximation
$$
\varepsilon_{\text{abs}} := | {\mathbb E}[Q({\bf y})] - {\mathcal
  A}|,
$$
where the estimator ${\mathcal A}$ is either ${\mathcal A}_{MFNNMC}$
or ${\mathcal A}_{HFMC}$. 
Given a $1\%$ failure probability ($\alpha=
0.01$) and a decreasing sequence of tolerances
$\varepsilon_{\footnotesize{\text{TOL}}} = 10^{-1}, 10^{-2}, 10^{-3}$,
a simple error analysis similar to the analysis in
Section \ref{sec:error} and verified by numerical computations gives
the minimum number of realizations $N$ and the maximum grid length
$h_{HF}$ for the high-fidelity model required to achieve $P( \varepsilon_{\text{abs}}  \le
\varepsilon_{\footnotesize{\text{TOL}}} ) = 0.99$. 
Table \ref{Ex2_table1} summarizes the numerical parameters $(N, h_{HF}, h_{LF})$ and the CPU time
of evaluating single realizations of $q_{LF}$ and $q_{HF}$.      
\begin{table}[!h]
        \centering
        \caption{Required number of realizations and grid lengths to
        achieve $P( \varepsilon_{\text{abs}}  \le
\varepsilon_{\footnotesize{\text{TOL}}} ) = 0.99$.}
\medskip
        \label{Ex2_table1}
        \begin{tabular}{|c|c|c|c|c|c|}
\hline
            $\varepsilon_{\footnotesize{\text{TOL}}}$ & $N$ & $h_{HF}$ &
          $W_{HF}$ & $h_{LF}$ & $W_{LF}$ \\
            \hline
            \hline
             $10^{-1}$ & $1.5 \times 10^{2}$ & $1/32$     & $0.67$ & $1/20$ & $0.21$\\
             $10^{-2}$ & $1.5 \times 10^{4}$ & $1/128$ & $29.75$ & $1/32$ & $0.67$\\
             $10^{-3}$ & $1.5 \times 10^{6}$ & $1/320$   & $708.21$ & $1/32$ & $0.67$\\
            \hline
        \end{tabular}
\end{table}

Following the algorithm in Section \ref{sec:algorithm}, we first
generate a uniform grid of $M = M_1 + M_2$ points ${\bf
  y}^{(i)} \in [10,11] \times [4,6]$, with
$i=1, \dotsc, M$, collected into two disjoint sets $Y_I$ and
$Y_{II}$. We select the two disjoint sets so that $M_2 \approx 3 \, M_1$, meaning that we
will need to compute the quantity $q({\bf y})$ by the high-fidelity
model with grid length $h_{HF}$ at only a quarter of points
$M/4$. The number of points $M$ will be chosen based on the desired
tolerance, slightly increasing as the tolerance decreases.

We will use the same architecture for the two networks NN1 and NN2 and
keep them fixed at all tolerance levels. Precisely, we choose feed-forward networks with
4 hidden layers, where each layer contains 30 neurons. We use ReLU
activation function for the hidden layers and the identity activation
function for the output layer of both networks. It is to be noted that
NN1 has three input neurons, while NN2 has two input neurons. Both networks
have one output neuron. For the training process, we split the
available $M$ data points into a training set (90$\%$ of $M$) and a
validation set (10$\%$ of $M$). 
We apply pre-processing transformations to the input data points
before they are presented to the two networks. Precisely, we transform the
points from $[10,11] \times [4,6]$ to the unit square $[0,1]^2$. 
We then employ the
quadratic cost function and use the Adam optimization technique with
an initial learning rate $\eta = 0.005$ that will be adaptively tuned
using the validation set. We do not use any regularization technique. Table \ref{Ex2_table2}
summarizes the number of training and validation data $M=M_1 + M_2$, the number of epochs $N_{\text{epoch}}$, batch size
$N_{\text{batch}}$, and the CPU time of training and evaluating the two
networks for different tolerances. 
We note that the number of training data satisfies $M
\propto \varepsilon_{\footnotesize{\text{TOL}}}^{-p}$ with $p=0.2$. 
\begin{table}[!h]
       \centering
        \caption{The number of training data and training and
          evaluation time of the two networks.}
\medskip
        \label{Ex2_table2}
\resizebox{\textwidth}{!}{
\begin{tabular}{|c||c|c||c|c|c|c||c|c|c|c|}
\hline
& & &\multicolumn{4}{c||}{NN1}&\multicolumn{4}{c|}{NN2}\\
\cline{4-11}
$\varepsilon_{\footnotesize{\text{TOL}}}$ & $M_1$ & $M_2$ &
                                                            $N_{\text{epoch}}$
                                          & $N_{\text{batch}}$ &
                                                                 $W_{T_1}$
                                          & $W_{P_1}$ &
                                                            $N_{\text{epoch}}$
                                          & $N_{\text{batch}}$ &
                                                                 $W_{T_2}$
                                          & $W_{P_1}$\\
\hline\hline
$10^{-1}$ & 848   & 2407   &500& 50& 326.57 & $4.10 \times
                                                          10^{-4}$ &
                                                                     500&
                                                                           50&
                                                                               1069.26 & $4.00 \times 10^{-4}$   \\
$10^{-2}$ & 1281 & 3680   &500& 50& 487.28 & $4.75 \times
                                                          10^{-4}$ & 1000& 50& 2904.67 & $4.60 \times 10^{-4}$   \\
$10^{-3}$ & 1976 & 5725 &1000& 50 & 1415.44 & $4.80 \times
                                                          10^{-4}$ &2000& 50& 9743.86 & $4.70 \times 10^{-4}$  \\
\hline
\end{tabular}}
\end{table}

\begin{figure}[!h]
\vspace{-.2cm}
\center
\subfigure{\includegraphics[width=8.1cm]{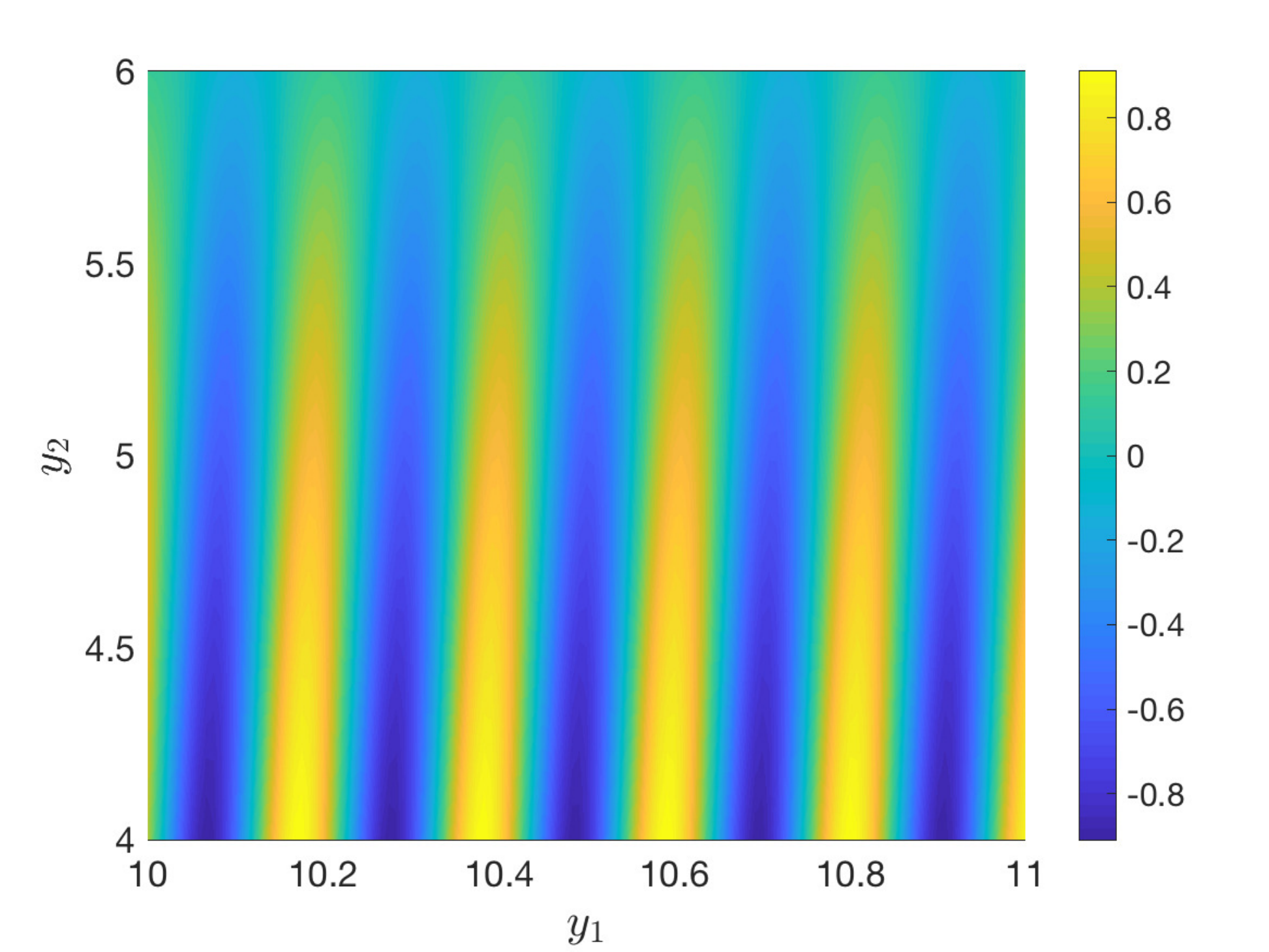}}
 \hskip .1cm
\subfigure{\includegraphics[width=8.1cm]{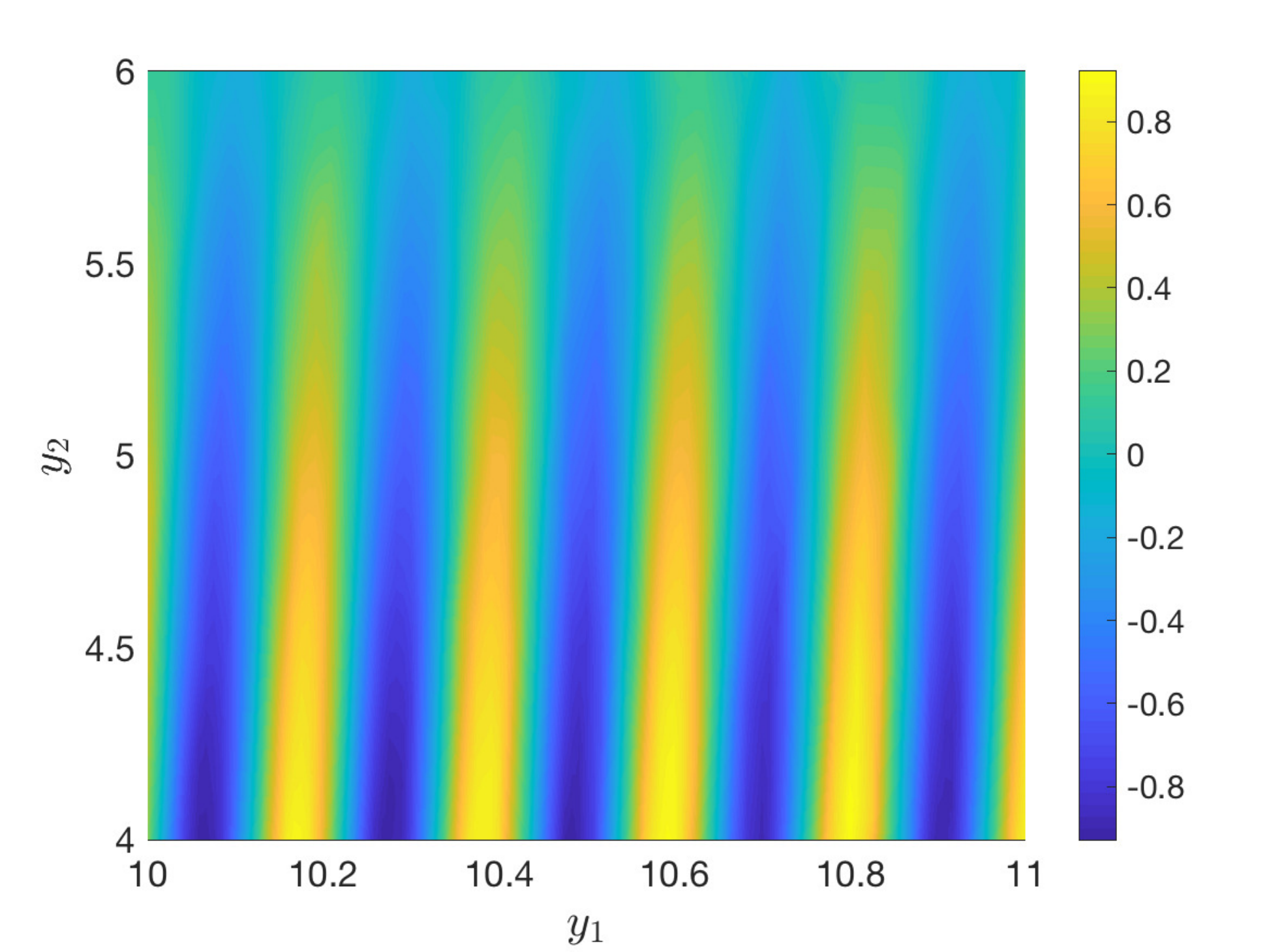}}
\vspace{-.2cm}
    \caption{High-fidelity quantity $q_{HF}({\bf y})$ for $\varepsilon_{\footnotesize{\text{TOL}}} =
  10^{-2}$. Left: true quantity. Right:
      predicted quantity by the trained network NN2.}
    \label{Ex2_toal2_predictions}
 \end{figure}

Figure \ref{Ex2_toal2_predictions} shows the true high-fidelity
quantity (left) and the predicted high-fidelity quantity by the
trained network NN2 (right) for tolerance $\varepsilon_{\footnotesize{\text{TOL}}}=10^{-2}$. 

\begin{figure}[!h]
\vspace{-0.2cm}
\begin{center}
\includegraphics[width=0.55\linewidth]{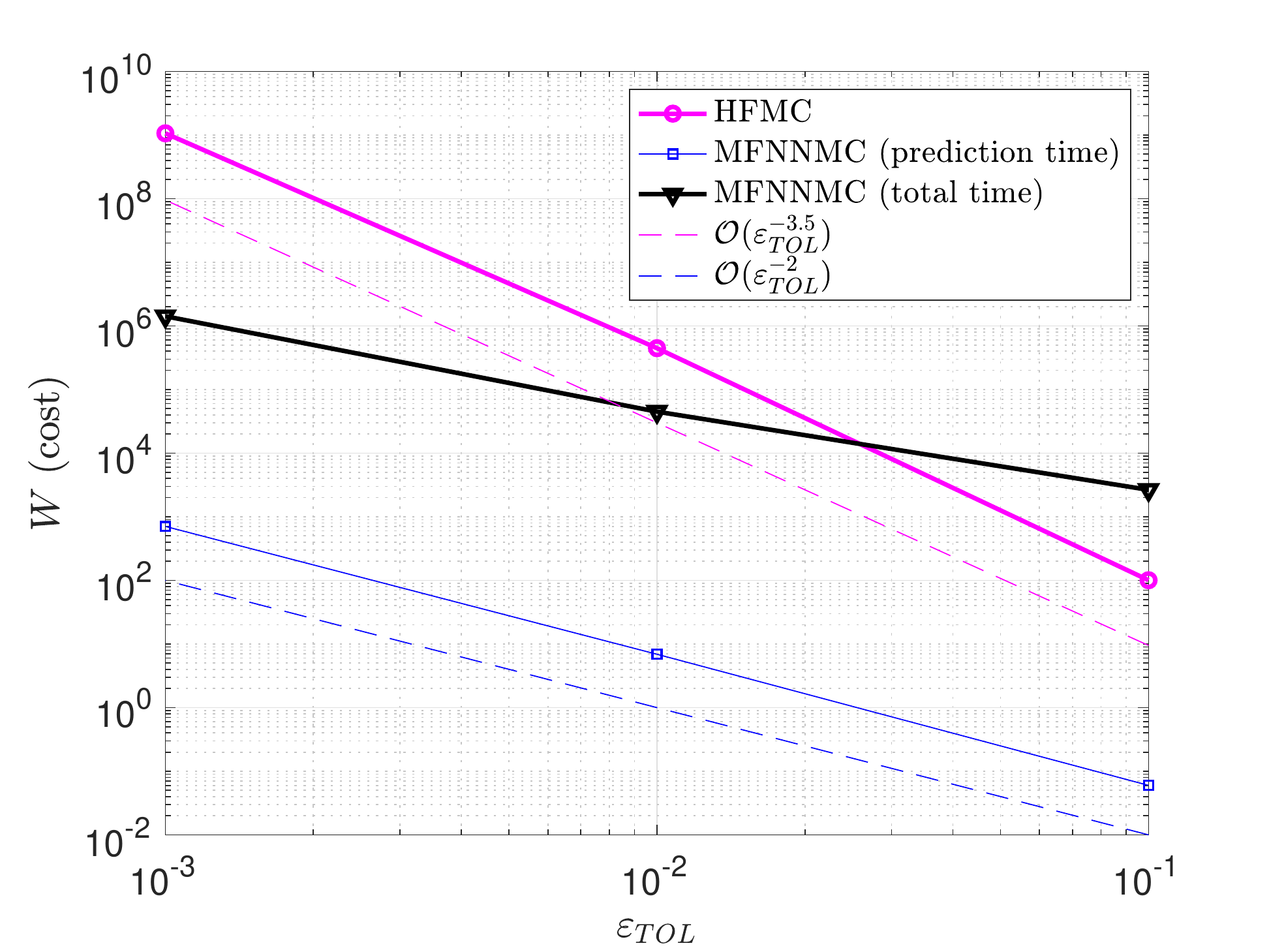}     
\vspace{-0.5cm}
\caption{CPU time versus tolerance. For large tolerances the
  training cost is dominant, making the cost of
  MFNNMC more than the cost of HFMC. However, as tolerance decreases,
  the training cost becomes negligible and the cost of
  MFNNMC approaches ${\mathcal
  O}(\varepsilon_{\footnotesize{\text{TOL}}}^{-2})$.}
\label{Ex2_conv}
\vspace{-.5cm}
\end{center}
\end{figure}

Figure \ref{Ex2_conv} shows the CPU time as a function of
tolerance. The computational cost of classical MC sampling is
proportional to $\varepsilon_{\footnotesize{\text{TOL}}}^{-3.5}$,
following \eqref{HFMC_cost_tol} and noting that the order of
accuracy of the finite difference scheme is $q=2$ and the time-space dimension of the problem
is $\gamma = 3$. 
On the other hand, if we only consider the prediction time of the
proposed multi-fidelity method, excluding the training costs, the cost of the proposed method is proportional to $\varepsilon_{\footnotesize{\text{TOL}}}^{-2}$ which
is much less than the cost of MC sampling. When adding the training
costs, we observe that although for large tolerances the training cost
is large, as the tolerance decreases the training costs become
negligible compared to the total CPU time. Overall, the cost of the
proposed method approaches ${\mathcal
  O}(\varepsilon_{\footnotesize{\text{TOL}}}^{-2})$ as tolerance
decreases, indicating orders of magnitude acceleration in computing
the expectation compared to MC sampling. This convergence rate can also be seen by \eqref{MFNNMC_cost_tol} where
$\max(2, p + \gamma/q) = \max(2, 0.2 + 1.5) = 2$. 

\begin{figure}[!h]
\begin{center}
\includegraphics[width=0.55\linewidth]{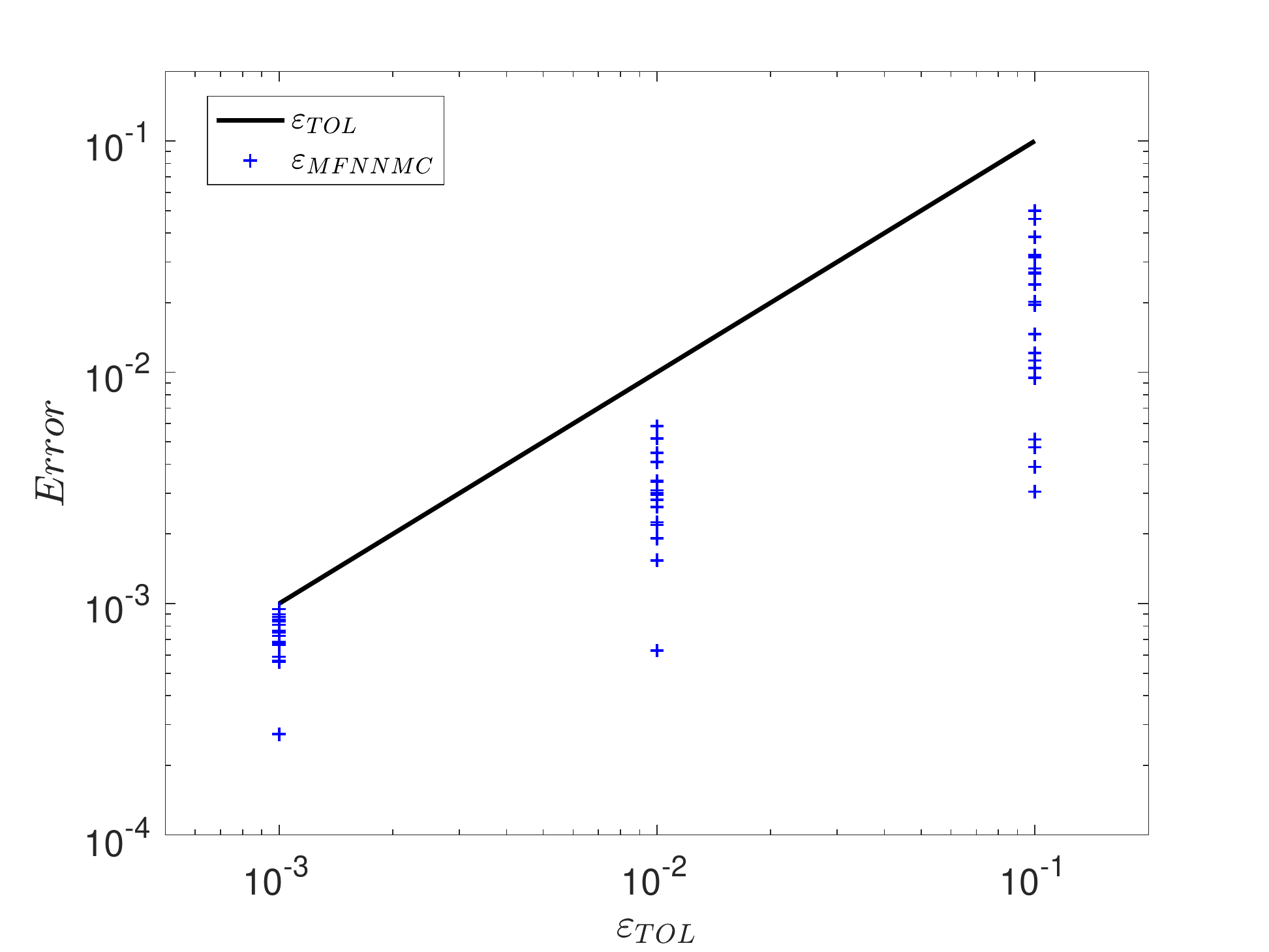}     
\vspace{-0.5cm}
\caption{Relative error as a function of tolerance, verifying that the
  tolerance is met with 1$\%$ failure probability. The ``+'' markers
  correspond to 20 simulations at each tolerance level.}
\label{Ex2_tols}
\end{center}
\end{figure}

Finally, Figure \ref{Ex2_tols} shows the relative error as a function
of tolerance for the proposed method, verifying that the tolerance
is met with 1$\%$ failure probability.

\section{Conclusion}
\label{sec:conclusion}

This work presents a multi-fidelity neural network surrogate sampling method
for the uncertainty quantification of physical/biological systems
described by systems of ODEs/PDEs. 
The proposed algorithm combines the approximation power of neural
networks with the advantages of MC sampling in a multi-fidelity
framework. 
For the numerical examples considered here, we observe
dramatic savings in computational cost when the output predictions are
desired to be accurate within small tolerances. 
More sophisticated numerical examples and a more comprehensive comparison
between the proposed method and other advanced MC sampling techniques
are subjects of current work and will be presented elsewhere. 
Other future directions include the extension of the proposed
construction  to training more than
two networks using data sets at multiple levels of fidelity within multi-level and multi-index frameworks.

\bibliographystyle{plain}
\bibliography{refs}

\end{document}